\documentclass[12pt]{article}
\usepackage{amsmath, amsthm,
amssymb}
\usepackage[all]{xy}
\title{Classification problems for system of forms and linear mappings%
\footnotetext{This is the authors' version of an article that was published in [\emph{Izv. Akad. Nauk SSSR Ser. Mat.} 51 (no. 6) (1987) 1170--1190] in Russian  and was translated into English in [\emph{Math. USSR-Izv.} 31 (no. 3) (1988) 481--501]. Theorem 2 in the article was formulated incorrectly. I correct it and add commentaries in footnotes.}
}
\author{Vladimir V. Sergeichuk\\ Institute of
Mathematics,
Tereshchenkivska 3,
Kiev,
Ukraine\\sergeich@imath.kiev.ua}

\date{}
\begin{document}

\renewcommand{\le}{\leqslant}
\renewcommand{\ge}{\geqslant}

\newcommand{\is}{\stackrel
{\text{\raisebox{-1ex}{$\sim\
\;$}}}{\to}}

\newcommand{\sdotss}%
{\text{\raisebox{-1.8pt}{$\cdot\,$}%
 \raisebox{1.5pt}{$\cdot$}%
\raisebox{4.8pt}{$\,\cdot$}}}

\newcommand{\sdotsss}%
{\text{\raisebox{-2.2pt}{$\cdot\,$}%
 \raisebox{1.7pt}{$\cdot$}%
\raisebox{5.6pt}{$\,\cdot$}}}

\newcommand{\diag}{\,\diagdown\,}
\newcommand{\dia}%
{\mathop{\rm diag}\nolimits}
\newcommand{\ind}%
{\mathop{\rm ind}\nolimits}

\newcommand{\End}%
{\mathop{\rm End}\nolimits}
\newcommand{\Imag}%
{\mathop{\rm Im}\nolimits}
\newcommand{\Ker}%
{\mathop{\rm Ker}\nolimits}
\newcommand{\Hom}%
{\mathop{\rm Hom}\nolimits}
\newcommand{\Aut}%
{\mathop{\rm Aut}\nolimits}
\newcommand{\im}{\mathop{\rm Im}\nolimits}
\newcommand{\rad}{\mathop{\rm rad}\nolimits}

\newcommand{\diagg}%
{\mathop{\rm diag}\nolimits}

\newtheorem{theorem}{Theorem}
\newtheorem{lemma}{Lemma}

\maketitle

\begin{abstract}
A method is proposed that allow the reduction of many classification problems of linear algebra to the problem of classifying Hermitian forms. Over the complex, real, and rational numbers classifications are obtained for bilinear forms, pairs of quadratic forms, isometric operators, and selfadjoint operators.
\end{abstract}

Many problems of linear algebra
can be formulated as problems of classifying the
representations of a quiver. A \emph{quiver} is, by definition, a directed graph.
A \emph{representation} of the quiver is given (see \cite{gab1}, and also \cite{b-gel-pon,naz}) by assigning to each vertex a
vector space and to each
arrow a linear mapping of the
corresponding vector spaces.
For example, the quivers
$$
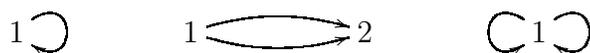
\xymatrix{
 {1}\ar@(ur,dr)
 &&
 {1}\ar@/^/[rr] \ar@/_/[rr] &&{2}
 &&
 {1} \ar@(ul,dl) \ar@(ur,dr)
 }
$$
correspond respectively to
the problems of classifying:
\begin{itemize}
  \item
linear operators (whose
solution is the Jordan or
Frobenius normal form),

  \item
pairs of linear mappings from
one space to another (the
matrix pencil problem, solved
by Kronecker), and
  \item
pairs of linear operators on
a vector space (a classical unsolved problem).
\end{itemize}
The notion of quiver has become central in the theory of finite-dimensional algebras over a field: the modules over an algebra are in one-to-one correspondence with the representations of a certain quiver with relations---the Gabriel quiver of the algebra (see \cite{gab1,rin}). The theories of quadratic and Hermitian forms are well developed (see \cite{mea,w.schar}).

We study systems of sesquilinear forms and linear mappings, regarding them as representations of a \emph{partially directed graph} (assigning to a vertex a vector space, to an undirected edge a sesquilinear form, and to a directed edge a linear mapping); and we show that the problem of classifying such representations, over a skew field $K$ of
characteristic $\ne 2$ reduces to the problems of classifying:

$1^{\circ}$ Hermitian forms over certain skew fields that are extensions of the center of $K$ and

$2^{\circ}$ representations of a certain quiver.

The quiver representations $2^{\circ}$ are known in the case of the problems of classifying:

(i)	\emph{bilinear} or \emph{sesquilinear forms} (see, for example, \cite{gab2,rie,rie1,sch}),

(ii)	\emph{pairs of symmetric}, or \emph{skew symmetric}, or \emph{ Hermitian forms} (\cite{r.schar1,w.schar,uhlig1,uhlig2,water,Yag}), and

(iii) \emph{isometric} or \emph{selfadjoint} operators on a space	 with nondegenerate symmetric, or skew symmetric, or Hermitian form (\cite{huppert1,huppert2,le,miln,sch,w.schar}).

We solve problems (i)--(iii) over a field $K$ of characteristic not $2$  with involution (possibly the identity) up to classification of Hermitian forms over fields that are finite extensions of $K$. This yields a \emph{classification of bilinear forms and pairs of quadratic forms over the rationals} since over finite extensions of the rationals classifications have been established for quadratic and Hermitian forms (see \cite{mea,w.schar}).

We study systems of forms and linear mappings by associating with them selfadjoint representations of a category with involution. This method was suggested by Gabriel \cite{gab2} for bilinear forms, and by Roiter \cite{roi} for systems of forms and linear mappings  (see also \cite{kru,ser_first}). Another approach to classification problems is proposed in \cite{q-s-s,w.schar}, where quadratic and Herinitian forms are studied on objects of an additive category with involution.

The main results of this paper were previously announced in \cite{ser_disch,ser_dep}.

The author wishes to thank A. V. Roiter for his considerable interest and assistance.

\section{Selfadjoint
representations of a linear
category with involution}
\label{s_category}

In this section we prove what might he called a weak Krull-Schmidt theorem for selfadjoint representations of a linear category with involution. \emph{Vector spaces are assumed throughout to be right vector spaces}.

By a \emph{linear category}
over a field $P$ is
meant a category $\cal C$ in
which for every pair of
objects $u, v$ the set of
morphisms $\Hom(u, v)$ is a
vector space over $P$
and multiplication of
morphisms is bilinear. The
set of objects in $\cal C$
will be denoted by ${\cal
C}_0$, the set of morphisms
by ${\cal C}_1$. We define
the category $R({\cal C})$ of
representations of ${\cal C}$
over a skew field $K$
with center $P$ as follows. A
\emph{representation} is a
functor $A$ from the category
${\cal C}$ to the category
${\cal V}$ of
finite-dimensional vector
spaces over $K$
having finite
dimension
\[
\dim(A) := \sum_{u\in{\cal
C}_0}\dim(A_u) < \infty
\]
and preserving linear
combinations:
\[
A_{\alpha a+\beta b}
=A_{\alpha} a+A_{\beta}
b,\qquad \alpha,\beta\in{\cal
C}_1,\ a,b\in P.
\]
(The images of an object $u$ and an morphism $\alpha$ are denoted by $A_u$ and $A_{\alpha}$.)
A \emph{morphism of
representations} $f\colon
A\to B$ is a natural
transformation of functors,
i.e., a set of linear
mappings
\[
f_u\colon A_u\to B_u,\qquad
u\in{\cal C}_0,
\]
such that
\[
f_uA_{\alpha}=B_{\alpha}f_u,\qquad
\alpha\colon u\to v.
\]

Suppose now that $K$
has an involution $a\mapsto
\bar a$; i.e., a bijection $K\to K$ satisfying $$\bar{\bar{a}}=a,\ \  \overline{a+b}= \bar{a}+\bar{b},\ \  \overline{ab}=\bar{b}\bar{a};$$ the involution can be the the identity if $K$ is a field. Following
\cite{roi}, we \emph{define
an involution on each of the
categories $\cal C$, $\cal
V$, and $R({\cal C})$}:
\begin{itemize}
  \item[1.]
To each object $u\in{\cal
C}_0$ we associate an object
$u^*\in{\cal C}_0$, and to
each morphism $\alpha\colon
u\to v$ a morphism
$\alpha^*\colon v^*\to u^*$
so that
\[
u^{**}=u\ne u^*,\qquad
\alpha^{**}=\alpha,\]\[
(\alpha\beta)^*=\beta^*\alpha^*,
\qquad (\alpha
a)^*=\alpha^*\bar a
\]
for all $u\in{\cal C}_0$,
$\alpha,\beta\in{\cal C}_1$,
$a\in P$ (note that
\cite{kru,roi} allow $u^* =
u$).

  \item[2.]
To each space $V\in\cal V$ we
associate the \emph{adjoint
space} $V^*\in\cal V$ of all
semilinear forms
$\varphi\colon V\to K$:
\[
\varphi(x+y)=\varphi(x)+
\varphi(y),\qquad
\varphi(xa)=\bar a\varphi(x)
\]
($x, y\in V;\ a\in K$), and to each linear
mapping $A: U\to V$ the
\emph{adjoint linear mapping}
$$A^*\colon V^*\to U^*,\qquad
A^*\varphi:=\varphi A.$$ We
identify $V$ and $V^{**}$.

  \item[3.]
To each representation $A\in
R({\cal C})$ we associate the
\emph{adjoint representation}
$A^{\circ}\in R({\cal C})$,
where
\[
A_u^{\circ}= A_{u^*}^*,\quad
A_{\alpha}^{\circ}=
A_{\alpha^*}^*\qquad
(u\in{\cal C}_0,\
\alpha\in{\cal C}_1);
\]
and to each morphism $f\colon
A\to B$ the \emph{adjoint
morphism} $f^{\circ}\colon
B^{\circ}\to A^{\circ}$,
where $f_u^{\circ}=
f_{u^*}^*$ ($u\in{\cal
C}_u$). An isomorphism
$f\colon A\to B$ of
selfadjoint representations
is called a \emph{congruence}
if $f^{\circ}=f^{-1}$.
\end{itemize}

In Section 2 we show that the problems of classifying the systems of sesquilinear forms and linear mappings over a skew field $K$ that satisfy certain relations with coefficients in the center of $K$ can be formulated as problems of classifying selfadjoint representations up to congruence. For the present we limit ourselves to examples.
\medskip

\noindent{\it Example 1.}
\[
{\cal C}_0=\{u,
u^*\},\qquad{\cal
C}_1=1_u{P}\cup
1_{u^*}{P}\cup
(\alpha{P}\oplus
\alpha^*{P}),
\]
where $\alpha,\alpha^*\colon
u\to u^*$. A selfadjoint
representation is given by a
pair of adjoint linear
mappings $A,A^*\colon U\to
U^*$, assigned to the
morphisms $\alpha$ and
$\alpha^*$. The
representation determines, in
a one-to-one manner, a
sesquilinear form $$A(x, y) :=
A(y)(x)$$ on the space $U$;
congruent representations
determine equivalent forms.
\medskip

\noindent{\it Example 2.}
\[
{\cal C}_0=\{u, u^*\},\qquad
{\cal C}_1=1_u{P}\cup
1_{u^*}{P}\cup
\alpha{P},
\]
where
$$\alpha=\varepsilon\alpha^*\colon
u\to u^*,\qquad
0\ne\varepsilon\in P.$$ A selfadjoint
representation determines an
$\varepsilon$-Hermitian form
$$A(x, y) = \varepsilon
\overline{A(y, x)}.$$

We show now how to obtain a
classification, up to
congruence, of the
selfadjoint representations
of a category $\cal C$,
starting with the knowledge
of a complete system
$\ind({\cal C})$ of its
nonisomorphic
direct-sum-indecomposable
representations. To begin
with, let us replace each
representation in $\ind({\cal
C})$ that is
\emph{isomorphic} to a
selfadjoint representation by
one that is \emph{actually}
selfadjoint, and denote the
set of such by $\ind_0({\cal
C})$. Denote by $\ind_1({\cal
C})$ the set consisting of
all representations in
$\ind({\cal C})$ that are
isomorphic to their adjoints
(but not to a selfadjoint),
together with one
representation from each pair
$ \{A,B\}\subset\ind({\cal
C})$ such that $A$ is not
isomorphic to $A^{\circ}$ but
is isomorphic to $B^{\circ}$.

In addition, we divide the
set ${\cal C}_0$ into two
disjoint subsets $S_0$ and
$S_0^*$ such that each pair
of adjoint objects $u, u^*$
has one member in $S_0$, the
other in $S_0^*$.

By the \emph{orthogonal sum}
$A\,\bot\, B$ of two
selfadjoint representations
$A$ and $B$ we mean the
selfadjoint representation
obtained from $A\oplus B$ by
specifying for each $v\in
S_0$ the action of
\[
\varphi+\psi\in A^*_v\oplus
B^*_v=(A\oplus B)_{v^*}
\]
on
\[
a+b\in A_v\oplus B_v
=(A\oplus B)_v
\]
as fol1ows:
\[
(\varphi+\psi)(a+b)=\varphi(a)
+\psi(b).
\]

For any representation $A$ we
define a \emph{selfadjoint
representation} $A^+$,
obtained from $A\oplus
A^{\circ}$ by specifying in a
similar fashion the action of

\begin{equation}\label{msf}
 (A\oplus A^{\circ})_{v^*}=
A_{v^*}\oplus A_v^*=
A_{v}^*\oplus A_{v^*}^{**}
\end{equation}
on
\[
(A\oplus A^{\circ})_{v}=
A_{v}\oplus A_{v^*}^*,\qquad
v\in S_0.
\]
Taking into account the
interchange of summands in
\eqref{msf}, we have
\begin{equation*}\label{ser01}
A_{\alpha}^+=\begin{bmatrix}
  0 &A_{\alpha}^{\circ} \\
A_{\alpha} & 0
\end{bmatrix}\qquad\text{for
}\ \alpha\colon u\to v^*,
\end{equation*}
\begin{equation*}\label{ser01a}
A_{\beta}^+=\begin{bmatrix}
  A_{\beta}&0\\
0& A_{\beta}^{\circ}
\end{bmatrix}\qquad\text{for
}\ \beta\colon u\to v,
\end{equation*}
\begin{equation*}\label{ser01b}
A_{\gamma}^+=\begin{bmatrix}
  0 &A_{\gamma} \\
A_{\gamma}^{\circ} & 0
\end{bmatrix}\qquad\text{for
}\ \gamma\colon u^*\to v,
\end{equation*}
where $u,v\in S_0$.

 For any selfadjoint
representation $A =
A^{\circ}$ and selfadjoint
automorphism $f = f^{\circ}$
of $A$, we define a
\emph{selfadjoint
representation $A^{f}$} and
an isomorphism
\begin{equation}\label{ser02a}
    \tilde f\colon A^f\to A,\qquad
\tilde f \tilde f^{\circ}=f,
\end{equation}
by putting
\[
 \tilde f_v=f_v,\quad
\tilde f_{v^*}=1\qquad (v\in
S_0)
\]
and
\[
 A^f_v=A_v,\quad
A_{\alpha}^f=\tilde
f_{v}^{-1}A_{\alpha}\tilde
f_{u}=1\qquad (v\in {\cal
C}_0,\ \ \alpha\colon u\to v)
\]

Now suppose $K$ has
characteristic $\ne 2$. We
show in Lemma \ref{lemser1a}
that the set $R$ of
noninvertible elements of the
endomorphism ring $$\Lambda =
\End(B),\qquad B \in
\ind_0({\cal C}),$$ is the
radical of $\Lambda$.
Therefore $T(B) = \Lambda/R$
is a skew field with
involution $$(f + R)^{\circ}
= f^{\circ}+ R.$$ For each
element $0\ne a =
a^{\circ}\in T(B)$ choose a
fixed automorphism
$f_a=f_a^{\circ}\in a$ (we
can take $f_a=(f +
f^{\circ})/2$, where $f\in
a$), and define $B^a =
B^{f_a}$. The set of
representations $B^a$ we call
the \emph{orbit of the
representation $B$}. For any
Hermitian form
\[
\varphi(x)=x^{\circ}_1a_1x_1+\dots+
x^{\circ}_ra_rx_r,\qquad 0\ne
a_i=a_i^{\circ}\in T({\cal
N}),
\]
we put
\[
B^{\varphi(x)}:=
B^{a_1}\,\bot\,\cdots\,\bot\,
B^{a_r}.
\]

\begin{theorem}\label{THEOREM 1}
Over a field or skew field
$K$ of characteristic
$\ne 2$, every selfadjoint
representation of a linear
category $\cal C$ with
involution is congruent to an
orthogonal sum
\begin{equation}\label{ser03}
A_1^+\,\bot\,\cdots\,\bot\,
A_m^+\,\bot\,B_1^{\varphi_1(x)}
\,\bot\,\cdots\,\bot\,
B_n^{\varphi_n(x)},
\end{equation}
where $A_i\in\ind_1({\cal
C})$, $B_j\in\ind_0({\cal
C})$, and $B_j\ne B_{j'}$ for
$j\ne j'$. The sum is
uniquely determined by the
original representation up to
permutation of summands and
replacement of
$B_j^{\varphi_j(x)}$ by
$B_j^{\psi_j(x)}$, where
$\varphi_j(x)$ and
$\psi_j(x)$ are equivalent
Hermitian forms over the skew
field $T(B_j)$.
\end{theorem}

\noindent {\it Remark}. Theorem 1 in fact holds for any ordinary (i.e., nonlinear) category $\cal C$ with involution, so long as we understand by a representation a functor $A:{\cal C}\to {\cal V}$ that has finite dimension $\dim(A) = \sum\dim(A_u)$. The ring $K$ can be replaced by any finite-dimensional quasi-Frobenius algebra $F$ with involution over a field of characteristic $\ne 2$ (a representation assigns to an object a finitely generated  module over $F$). A finite-dimensional algebra $F$ is quasi-Frohenius if the regular module $F_F$ is injective; over such an algebra the finitely generated modules $ M$ and $M^{**}$ can still be identified.
\medskip

Theorem \ref{THEOREM 1}
reduces the classification,
up to congruence, of the
selfadjoint representations
of the category $\cal C$,
assuming known the
representations $\ind_1({\cal
C})$ and the orbits of the
representations $\ind_0({\cal
C})$, to the classification
of Hermitian forms over the
skew fields $$T(B),\qquad
B\in \ind_0({\cal C}).$$
\emph{If $K$ is a
finite-dimensional over its
center $Z$, then $T = T(B)$
is finite-dimensional over
$Z$ under the natural
imbedding of $Z$ in the
center of $T$, and the
involution on $T$ extends the
involution on $Z$.}

Suppose, for example, that
$K$ is a
\emph{real closed field};
i.e., $$1 < (K_{\text{alg}}:K)
<\infty,$$ where $K_{\text{alg}}$ is the
algebraic closure of $K$. Then its characteristic
is $0$, $K_{\text{alg}} =K(\sqrt{-1})$, and $K$ has only the identity
involution: the stationary
subfield relative to
involution must coincide with
$K$ (see
\cite[Chap. VI,
\S\,2, nos. 1, 6 and Exercise
22(d)]{bour}). By the theorem of
Frobenius \cite{bour2}, $T$
is equal to either $K$, or
$K_{\text{alg}}$, or
the algebra $\mathbb H$ of
quaternions over $K$.
By the law of inertia
\cite{bour2}, if $T = K$ or $T=K_{\text{alg}}$ with
nonidentity involution, or $T
= \mathbb H$ with the
\emph{standard involution}
\[
a+bi+cj+dk\mapsto a-bi-cj-dk,
\]
then
a Hermitian form over $T$ is
equivalent to exactly one
form of the form
\[
x_1^{\circ}x_1+\dots
+x_l^{\circ}x_l
-x_{l+1}^{\circ}x_{l+1}
-\dots- x_r^{\circ}x_r.
\]
If $T = \mathbb H$ with
nonstandard involution, then
every Hermitian form
\[
\varphi(x)=x_1^{\circ}a_1x_1+\dots
+x_r^{\circ}a_rx_r\qquad
(a_i=a_i^{\circ}\ne 0)
\]
over $T$ is equivalent to the
form
\[
x_1^{\circ}x_1+\dots
+x_r^{\circ}x_r
\]
since $$a_i = b_i^2 =
b_i^{\circ}b_i,$$ where
$b_i\in K(a_i)$ for
$a_i\notin K$ and
$b_i\in K(d)$ for
$a_i\in K$. Here
$d=d^{\circ}\notin K$;
the existence of $d$ follows
from \cite{bour2} (Chap.
VIII, \S\,11, Proposition 2);
and $K(a_i)$ and
$K(d)$ are
algebraically closed fields
with the identity involution.

Suppose $K$ is a
\emph{finite field}. Then $T$
is also a finite field, over
which a Hermitian form
reduces uniquely to the form
\[
x_1^{\circ}tx_1+
x_2^{\circ}x_2+\dots
+x_r^{\circ}x_r,
\]
where $t$ is equal to 1 for
nonidentity involution on
$T$, and $t$ is equal to 1 or
a fixed nonsquare for the
identity involution
(\cite[Chap. 1, \S\,8]{die}).

Thus, applying Theorem
\ref{THEOREM 1}, we obtain
the following assertion, a
special case of which is the
law of inertia for quadratic
and Hermitian forms.

\begin{theorem}\label{THEOREM 2}%
Let $K$ be one of the
following fields or skew
fields of characteristic not
$2$:
\begin{itemize}
  \item[\rm(a)]
An algebraically closed field
with the identity involution.

  \item[\rm(b)]
An algebraically closed field
with nonidentity involution.

  \item[\rm(c)]
A real closed field or
the algebra of quaternions
over a real closed field.

  \item[\rm(d)]
 A finite field.
\end{itemize}

Then over $K$ every
selfadjoint representation of
a linear category $\cal C$
with involution is congruent
to an orthogonal sum,
uniquely determined up to
permutation of summands, of
representations of the
following form $($where
$A\in\ind_1({\cal C})$ and
$B\in\ind_0({\cal C}))$:
\begin{itemize}
  \item[\rm(a)]
$A^+$, $B$.

  \item[\rm(b)]
$A^+$, $B$, $B^{-1}$
$(-1\in\Aut(B))$.

  \item[\rm(c)]
$A^+$, $B^{t}$, where $t = 1$
if $T(B)$ is an algebraically
closed field with the
identity involution or the
algebra of quaternions with
nonstandard involution, and
$t = ±1$ otherwise.

  \item[\rm(d)]
$A^+$, $B^{t}$, where $t$ is
equal to $1$ for nonidentity
involution on the field
$T(B)$, and is equal to $1$
or a fixed nonsquare in
$T(B)$ for the identity
involution, while for each
$B$ the orthogonal sum has at
most one summand $B^t$ with
$t\ne 1$.\footnote{This theorem was formulated in [V. V. Sergeichuk,
Classification
problems  for systems
of forms and linear
mappings, {\it  Math.
USSR-Izv.} 31 (no. 3)
(1988) 481--501] incorrectly in the case of the algebra $\mathbb
H$ of quaternions over a real closed field.
Formulating it, I
erroneously thought
that all $T(B)=\mathbb H$
if $K=\mathbb
H$. I wrote that the summands have the form (a) or (b) if $K=\mathbb H$ and the involution on $\mathbb H$ is nonstandard or standard, respectively.}
\end{itemize}
\end{theorem}

\noindent{\it Remark 1.}
A similar assertion can be
made for representations over
a field $K$ of
algebraic numbers since over
skew fields that are finite
central extensions of such a
$K$ the
classification of Hermitian
forms is known $($see
{\rm\cite[Chap. 1,
\S\,8]{die}}$)$.
\smallskip

\noindent{\it Remark 2.} It can be shown that over an algebraically closed field of characteristic $\ne 2$, or over a real closed field, a system of tensors of valence $\ge 2$ decomposes uniquely,
to isomorphism of summands, into a direct sum of indecomposable subsystems. For a system of valence $2$ this hollows from Theorem 2 (see Section 2) Over an algebraically closed field of characteristic 2, on the other hand, even the number of summands depends on particular decomposition: the symmetric bilinear forms $x_1y_1 + x_2y_2 +x_3y_3$ and
$x_1y_2 + x_2y_1 +x_3y_3$ are equivalent, hut the form $x_1y_2 + x_2y_1$ is indecomposable.
\medskip

The rest of this section is
devoted to the proof of
Theorem \ref{THEOREM 1}.

\begin{lemma}\label{lemser1a}
The radical of the
endomorphism ring of a
direct-sum-indecomposable
representation is a nilpotent
ideal and consists of all
noninvertible endomorphisms.
\end{lemma}

\begin{proof}
For every
representation $B$, any
endomorphism $f$ of $B$
satisfies Fitting's lemma:
$${B} = \Imag(f^d)\oplus
\Ker(f^d),\qquad d =
\dim(B),$$ where
$\Imag(f^d)$ and $\Ker(f^d)$
are the restrictions of $B$ to the subspaces
$\Imag(f_v^d)$ and
$\Ker(f_v^d)$ for all
$v\in{\cal C}_0$.

Suppose $B$ is
direct-sum-indecomposable,
and let $R$ be its set of
noninvertible endomorphisms.
Then $f^d = 0$ for all $f\in
R$. If $$\Imag(f) =
\Imag(fg)\qquad (f, g\in
R),$$ then $$\Imag(f) =
\Imag(fg^d) = 0,$$ i.e., $f =
0$. Consequently, $$f_1\cdots
f_d = 0,\qquad (f_1 + f_2)^d
= 0$$ for all $f_i\in R$, and
$R$ is a nilpotent ideal of
the endomorphism ring and
coincides with the radical.
\end{proof}

\begin{lemma}\label{lemser2a}
Let $A$ be a selfadjoint
representation of a category
$\cal C$ over a skew field
$K$ of characteristic
$\ne 2$ that is
indecomposable into an
orthogonal sum but possesses
a nontrivial decomposition
into a direct sum. Then $A$
is congruent to $B^+$, where
$B\in\ind_1({\cal C})$.
\end{lemma}

\begin{proof}
$1^{\circ}$. Let $f =\pm
f^{\circ}$ be an endomorphism
of $A$. By Fitting's lemma,
$$A = \Imag(f^d)\oplus \Ker(f^d),
\qquad d =\dim(A).$$ But $A$
is indecomposable into an
orthogonal sum. Therefore $f$
is either invertible or
nilpotent.
\medskip

$2^{\circ}$. Since $A$ is
direct-sum-decomposable,
there exists a nontrivial
idempotent $e = e^{\circ}\in
\End(A)$, the projection onto
an indecomposable direct
summand. From $1^{\circ}$ it
follows that $ee^{\circ}$ is
nilpotent. Consider the
selfadjoint endomorphism $h =
p(ee^{\circ})$, where
\begin{equation}\label{ser04}
p(x) = 1 + a_1x + a_2x^2
+\cdots
\end{equation}
is an infinite series with
coefficients in the prime
subfield of $K$ such
that $p(x)^2 =1 - x.$ The
series exists, since the
characteristic $\ne 2$.

Consider the idempotent $f=
h^{-1}eh$. It satisfies
\begin{align*}
ff^{\circ} &=
h^{-1}eh^2e^{\circ}h^{-1} =
h^{-1}e(l -
ee^{\circ})e^{\circ}h^{-1}
\\ &= h^{-1}(ee^{\circ} -
ee^{\circ})h^{-1} = 0.
\end{align*}
If we take a new $e$ equal to
$f^{\circ}$, then $e^{\circ}e
= 0$ and for the new $f =
h^{-1}eh$ we find, in
addition to $ff^{\circ} = 0$,
that
\begin{align*}
f^{\circ}f &= he^{\circ}
h^{-2} eh = he^{\circ}(1 -
ee^{\circ})^{-1}eh \\ &=
he^{\circ}(1 + ee^{\circ})eh
= 0.
\end{align*}
Consequently, $f +f^{\circ}$
is idempotent. By
$1^{\circ}$, $f+f^{\circ}$ is
a trivial idempotent. If
$f+f^{\circ} = 0$, then
$f^{\circ} = -f$, a
contradiction to $1^{\circ}$.
Therefore $f + f^{\circ} = 1$
and $A = D^+$, where the
representation $D = \Imag(f)$
is direct-sum-indecomposable
(recall that $e$ is the
projection onto an
indecomposable direct
summand).

The representation $D$ cannot
be isomorphic to one that is
selfadjoint. Indeed, suppose
$$\varphi\colon D \to B =
B^{\circ}$$ is an
isomorphism. Then
$$\varphi\oplus
(\varphi^{\circ})^{-1}\colon
D^+\to B^+$$ and
$$\psi\colon B^+\to
B\,\bot\, B^{-1}\qquad (-1\in
\Aut(B))$$ are congruences,
where
\[
\psi_v=\begin{bmatrix}
  1/2 &1 \\
-1/2 & 1
\end{bmatrix},\qquad
\psi_{v^*}=\psi_v^{*-1}\qquad
(v\in S_0).
\]
This contradicts the
assumption that $A = D^+$ is
indecomposable into an
orthogonal sum.

Consequently, the
representation $D$ is
isomorphic either to $B$ or
to $B^{\circ}$, where
$B\in\ind_1({\cal C})$. Then
$D^+$ is congruent either to
$B^+$ or to $(B^{\circ})^+$.
But $(B^{\circ})^+$ is
congruent to $B^+$. Therefore
$A = D^+$ is congruent to
$B^+$; and the lemma is
proved.
\end{proof}

\noindent{\it Remark.}
Over a skew field of
characteristic 2 the lemma is
false, but a weaker version
holds: if $A$ is the
representation of the lemma,
then $A\simeq B \oplus
B^{\circ}$, where $B$ is
indecomposable. Indeed, let
\[
g = e + e^{\circ} +
ee^{\circ}\in \End(A),
\]
where $e = e^2$ is the
projection onto the direct
summand $B$ of least
dimension $\sum_u \dim(B_u)$.
If $g$ is noninvertible, it
is nilpotent (step
$1^{\circ}$ in the proof of
the lemma) and $h = 1 + g$ is
invertible; which implies,
since
\[
eh = e + e + ee^{\circ} +
ee^{\circ} = 0,
\]
that $e = 0$. Therefore $g$
must be invertible, and
\[
A =(e+e^{\circ}+ee^{\circ})A
=eA+e^{\circ}A\simeq B\oplus
B^{\circ}.
\]

\begin{lemma}\label{lemser3a}
Let $A$ and $B$ be two
selfadjoint representations,
\[
f = f^{\circ}\in
\Aut(A),\qquad g =
g^{\circ}\in\Aut(B).
\]
Then $A^f$ is congruent to
$B^g$ if and only if $g =
h^{\circ}fh$ for some
isomorphism $h\colon B\to A$.
\end{lemma}

\begin{proof}
Suppose $\varphi\colon B^g\to
A^f$ is a congruence. Define
the isomorphism
\[
h = (\tilde
f^{\circ})^{-1}\varphi \tilde
g^{\circ}\colon B\to A,
\]
where $\tilde f$ and $\tilde
g$ are the isomorphisms of
form \eqref{ser02a}. Using
the relations \[f = \tilde
f\tilde f^{\circ},\qquad g =
\tilde g\tilde
g^{\circ},\qquad
\varphi^{\circ}\varphi= 1,
\]
we find $h^{\circ}fh = g$.

Conversely, if $g =
h^{\circ}fh$, then
\[
\varphi = \tilde
f^{\circ}h(\tilde
g^{\circ})^{-1}\colon B^g\to
A^f\] is a congruence.
\end{proof}

\begin{lemma}\label{lemser4a}
Suppose a representation $A$
is selfadjoint and
direct-sum-indecomposable.
Then $A$ is congruent to a
representation $B^f$, where
$B\in\ind_0({\cal C})$ and $f
= f^{\circ}\in \Aut(B)$.
\end{lemma}

\begin{proof}
By definition of the set
$\ind_0({\cal C})$, there
exists an isomorphism
$h\colon B\to A$ for some
$B\in\ind_0({\cal C})$. By
Lemma \ref{lemser3a},
$B^{h^{\circ}h}$ is congruent
to $B$.
\end{proof}

\begin{lemma}\label{lemser5a}
Over a skew field $K$
of characteristic not $2$, the
representations
\[
B^{f_1}\,\bot\,\cdots\,\bot\,
B^{f_n}\text{\quad and \quad}
B^{g_1}\,\bot\,\cdots\,\bot\,
B^{g_n}
\]
\[(B\in \ind_0({\cal C}),\quad
f_i~ =
f_i^{\circ}\in\Aut(B),\quad
g_i =
g_i^{\circ}\in\Aut(B))\] are
congruent if and only if the
Hermitian forms
\[
\sum x_i^{\circ}(f_i +
R)x_i,\qquad \sum
x_i^{\circ}(g_i + R)x_i
\]
are equivalent over the skew
field $T(B) = \End(B)/R$.
\end{lemma}

\begin{proof}
Obviously,
\[
B^{f_1}\,\bot\,\cdots\,\bot\,
B^{f_n}= D^f,
\]
where
\[
D:=B\,\bot\,\cdots\,\bot\,
B,\qquad f:=
\dia(f_1,\dots,f_n).
\]
By Lemma \ref{lemser3a},
$D^f$ is congruent to $D^g$
for some $$g =
\dia(g_1,\dots, g_n)$$ if and
only if $g = h^{\circ}fh$,
where $h = [h_{ij}],\
h_{ij}\in \End(B)$.

In particular, if $g_i\in f_i
+ R$, then $D^g$ is congruent
to $D^f$. Indeed, let $h =
p(r)$, where $p(x)$ is the
series \eqref{ser04} and
\[
r =
\dia(r_1,\dots,r_n),\qquad
r_i= f_i^{-1}(f_i -g_i)\in R
\]
(by Lemma \ref{lemser1a}, the
matrix $r$ is nilpotent).
Then
\[
r_i^{\circ}f_i = (1 - g_i
f_i^{-1})f_i=f_ir_i,
\]
\[
h^{\circ}fh = p(r^{\circ})fh
= fh^2 = f(1 - r) =g.
\]

Consequently, all the
matrices in the set
$$\dia(f_1 + R,\dots , f_n +
R)$$ give congruent
representations. Thus, $D^g$
is congruent to $D^f$ if and
only if
\[
\dia(b_1,\dots, b_n)
=[c_{ij}]^{\circ}
\dia(a_1,\dots, a_n)
[c_{ij}],
\]
where
\[
a_i=f_i+R,\quad
b_i=g_i+R,\quad c_{ij}=h_{ij}
+ R;
\]
i.e., if and only if the
Hermitian forms
\[
\sum x_i^{\circ}a_ix_i,\qquad
\sum x_i^{\circ}b_ix_i
\]
are equivalent over $T(B)$.
\end{proof}

\begin{proof}[Proof of Theorem
{\ref{THEOREM 1}}]
$1^{\circ}$. If
\[
M_1\oplus\dots\oplus M_t
\overset{f}{\to}
N\oplus\dots\oplus N
\overset{g}{\to}
M_1\oplus\dots\oplus M_t
\]
are two homomorphisms of
direct sums of indecomposable
representations of the
category $\cal C$, with $N$
nonisomorphic to any of the
representations $M_1,\dots,
M_t$, then the endomorphism
$h:=fg$ is nilpotent. Indeed,
$f,\ g$, and $h$ can be
written as matrices:
\[
f=[f_{ij}],\qquad g =
[g_{jk}],\qquad h = [h_{ik}],
\]
where
\[
f_{ij}\colon M_j\to N,\qquad
g_{jk}\colon N \to M_j,\] and
\[h_{ik} =\sum_j
f_{ij}g_{jk}\colon N \to N.
\]
Since the set $R$ of
non-invertible elements of
the ring $\End(N)$ is a
nilpotent ideal (Lemma
\ref{lemser1a}), it suffices
to show that $f_{ij}g_{ik}\in
R$. Suppose that, on the
contrary, $f_{ij}g_{ik}$ is
invertible. Then so is
$g_{ik}f_{ij}$ (since it is
not nilpotent); and therefore
$f_{ij}$ is an isomorphism,
contradicting the assumption
$M_{j}\not\simeq N$.

$2^{\circ}$.  By Lemmas
\ref{lemser2a} and
\ref{lemser4a}, every
selfadjoint representation is
congruent to a representation
$A$ of the form
(\ref{ser03}). Let
\begin{equation}\label{ser05}
C=C_1^+\,\bot\,\cdots\,\bot\,
C_k^+\,\bot\,D_1^{\psi_1(x)}
\,\bot\,\cdots\,\bot\,
D_l^{\psi_l(x)}
\end{equation}
be another representation of
the same form, and $f\colon
A\to C$ a congruence. Since
the representations $A$ and
$C$ are isomorphic, so are
their indecomposable direct
summands (the Krull-Schmidt
theorem \cite[Chap. I,
Theorem (3.6)]{bas} for the
additive category $R({\cal
C}))$. In view of the
isomorphism
\[
A\simeq
\bigoplus_{i=1}^m(A_i\oplus
A_i^{\circ})\oplus
\bigoplus_{j=1}^n(B_j\oplus
\dots\oplus B_j)
\]
(see (\ref{ser03}) and
(\ref{ser02a}))  and the
definition of the sets
$\ind_0({\cal C})$ and
$\ind_1({\cal C})$, we find
that $m = k$ and $n = l$, and
that, reindexing if
necessary,
\[
A_i = C_i,\qquad B_j =
D_j,\qquad
B_j^{\varphi_j(x)}\simeq
B_j^{\psi_j(x)}.
\]

Write the congruence $f\colon
A \to C$ as a matrix
\[
f=
\begin{bmatrix}
  f_{11} &f_{12} \\
f_{21} & f_{22}
\end{bmatrix}\colon S\,\bot\,
B_n^{\varphi_n(x)}\to
T\,\bot\, B_n^{\psi_n(x)},
\]
where $S$ and $T$ are the
sums (\ref{ser03}) and
(\ref{ser05}) without the
last summand. From
$f^{\circ}f = 1$ it follows
that
\[
f_{12}^{\circ}f_{12} +
f_{22}^{\circ}f_{22} = 1.
\]
Since $f_{12}^{\circ}f_{12}$
is nilpotent (see
$1^{\circ}$), we can define
the homomorphism
\[
g =
f_{22}p(f_{12}^{\circ}f_{12})^{-1}
\colon B_n^{\varphi_n(x)}\to
B_n^{\psi_n(x)},
\]
where $p(x)$ is the series
(\ref{ser04}). Since
\[
g^{\circ}g =
p(f_{12}^{\circ}f_{12})^{-1}(1
- f_{12}^{\circ}f_{12})
p(f_{12}^{\circ}f_{12})^{-1}
= 1
\]
and
\[
B_n^{\varphi_n(x)}\simeq
B_n^{\psi_n(x)},
\]
$g$ is a congruence. By Lemma
\ref{lemser5a}, the Hermitian
forms $\varphi_n(x)$ and
$\psi_n(x)$ are equivalent. A
similar argument gives
equivalence of each of the
forms $\varphi_j(x)$ and
$\psi_j(x)$  $(1 \le j <n)$.
\end{proof}

\section{Applications to linear algebra}

In this section we apply Theorem 1 to some classical classification problems.

Let $\cal C$ be a linear category with involution over a field $P.$ To specify the category $\cal C$, it suffices to list:

\begin{itemize}
  \item[\rm(i)]
a set $S_0\in{\cal C}_0$ of
objects of the category such
that
\[
S_0\cup S_0^*={\cal
C}_0,\qquad S_0\cap
S_0^*=\varnothing,
\]

  \item[\rm(ii)]
a set $S_1\in{\cal C}_1$ of
generating morphisms, such
that every morphism in the
category is representable as
a linear combination of
products of morphisms in
\[
S_1\cup S_1^*\cup
\{1_u\,|\,u\in {\cal C}_0\},
\]
  \item[\rm(iii)]
a set $S_2$ of defining
relations for ${\cal C}$:
\[
\sum_i\alpha_{i1}\cdots
\alpha_{it_i} a_i = 0,
\]
\[
(a_i\in P,\qquad
\alpha_{ij}\in S_1\cup
S_1^*\cup \{1_u\,|\,u\in
{\cal C}_0\}),
\]
\end{itemize}
such that multiplication of
morphisms in ${\cal C}$ is
completely determined by the
bilinearity property and the
relations $S_2 \cup S_2^*$,
where $S_2^*$ consists of the
adjoints of the relations in
$S_2$:
\[
\sum_i\alpha_{it_i}^*\cdots
\alpha_{i1}^*\bar a_i = 0.
\]

Let us agree, further, that
\emph{the set $S_1$ does not
contain any morphisms of the
form $\alpha\colon v^*\to
u^*$ $(u, v \in
S_0)$}---since these can be
replaced by the adjoint
morphisms $\alpha^*\colon
u\to v$.

If the sets $S_0$ and $S_1$
are finite, $S_2$ are also be
taken to be finite. Such
categories, called finitely
generated, can be
conveniently presented by
graphs in the following two
ways:
\begin{itemize}
  \item
By a quiver $\overline S$
with the set of vertices
${\overline S}_0:= S_0\cup
S_0^*$, the set of arrows
${\overline S}_1:= S_1\cup
S_1^*$, and the set of
defining relations
${\overline S}_2:= S_2\cup
S_2^*$. Such a quiver is
called a \emph{quiver with involution}
of the category $\cal C$ (see
\cite{roi}).

  \item
By a graph
$S$
with the set of vertices
$S_0$, the set of edges
$S_1$, and the set of
defining relations $S_2$.
Each morphisms in $S_1$ of
the form
\[
\alpha\colon u \to v*,\quad
\beta\colon u \to v,\quad
\gamma\colon u^* \to v\quad
(u,v\in S_0)
\]
is represented, respectively,
by edges of the form
\[
\alpha\colon
u\:\frac{\quad}{}\:v,\quad
\beta\colon u \to v,\quad
\gamma\colon
u\leftrightarrow v.
\]
Such a graph, with nondirected, directed, and doubly directed edges, we call a \emph{doubly oriented graph} (\emph{dograph} for short) of the category $\cal
C$.%
\footnote{Such a graph is called a \emph{discheme} (directed
scheme) in \cite{ser_first} and in [V. V. Sergeichuk,
Classification
problems  for systems
of forms and linear
mappings, {\it  Math.
USSR-Izv.} 31 (no. 3)
(1988) 481--501].}
\end{itemize}

For example:
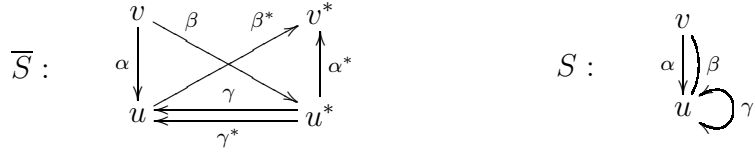
\[
\raisebox{23pt}{\xymatrix@R=4pt{
 &{v}\ar[dd]_{{\alpha}}
 \ar[ddrr]^(.25){{\beta}}&
 &{v^*} \\
 {\overline{S}:}&&\\
 &{u}\ar[uurr]^(.75){{\beta}^*}
 \ar@<0.4ex>[rr];[]^{{\gamma^*}}
 \ar@<-0.4ex>[rr];[]_{{\gamma}}
 &&{u^*}\ar[uu]_{{\alpha}^*}
 }}\qquad\qquad\qquad
\raisebox{20pt}{\xymatrix@R=4pt{
 &{v}\ar[dd]_{\alpha}
 \ar@{-}@/^/[dd]^{\beta}\\
{S}:&
  \\
 &{u} \ar@(ur,dr)@{<->}^{\gamma}}}
\]

In what follow, a representation of the
category ${\cal C}$ will be
specified not on the whole
set ${\cal C}_0\cup {\cal
C}_1$, but on the subset
${\overline S}_0\cup {\overline S}_1$ (being
completely determined by its
values on the subset); and we shall speak, correspondingly, not of a representation of the category $\cal C$ but of a representation of the quiver $\overline S$. Thus, a representation $A$ of the quiver $\overline S$ over a skew field $K$ is a set of finite-dimensional vector spaces $A_v$ $(v\in {\overline S}_0$) over $K$ and linear mappings $A_{\alpha}: A_u \to A_v$ $({\overline S}_1\ni \alpha:u\to v$) satisfying the relations ${\overline S}_2$ (with the $\alpha\in{\overline S}_1$ replaced by the $A_{\alpha}$).
A selfadjoint representation is completely determined by its values on the set $S_0\cup S_1$, i.e., by a set of finite-dimensional vector spaces $A_v$ ($v\in S_0$) and linear mappings $A_{\alpha}$ ($\alpha\in S_1$) of the form $A_u \to A_v^*$ for $\alpha : u \frac{\quad}{} v$, $A_u \to A_v$ for $\alpha : u \to v$, and $A_u^*\to A_v$ for $\alpha : u \leftrightarrow v$, satisfying the relations $S_2$ (with $\alpha \in S_1$ and $\alpha^* \in S_1^*$ replaced by $A_{\alpha }$ and $A_{\alpha }^*$). Such a set will he called a \emph{representation $A$ of the dograph $S$} (see \cite{ser_first}).

\emph{A linear mapping $A: U \to V^*$ will be identified with the sesquilinear form $A: V \times U \to K$, $A(v,u):=A(u)(v)$} (their matrices coincide, with the understanding that in the adjoint space we choose the adjoint basis). Recall that by a sesquilinear form is meant a mapping $A: V\times U\to K$ such that
\begin{gather*}
  A(va+v'a',u)=\bar{a}A(v,u)+
  \bar{a'}A(v',u),\\
  A(v,ua+u'a')=A(v,u)a+A(v,u')a'
\end{gather*}
for all $v,v'\in V$, $u,u'\in
U$, and $a,a'\in K$.

With this identification, a representation $A$ of a dograph $S$ is a set of vector spaces $A_v$ ($v\in S_0$), and linear mappings and sesquilinear forms $A_{\alpha}$ ($\alpha\in S_1$) of the form
$A_v \times A_u\to K$ for $\alpha : u \frac{\quad}{} v$, $A_u \to A_v$ for $\alpha : u \to v$, and $A_v^*\times A_u^*\to K$ for $\alpha : u \leftrightarrow v$ (in other words, the $A_{\alpha }$ are doubly covariant, mixed, or doubly contravariant tensors on the spaces $A_u$ and $A_v$). For two representations $A$ and $B$ of a dograph there is also a natural translation of the notions of
\begin{itemize}
  \item
\emph{congruence} $f: A\to B$---a set of nonsingular linear mappings $f_v: A_v \to B_v$ ($v \in S_0$) taking the $A_{\alpha }$, into the $B_{\alpha }$ ($\alpha \in S_1$)---and of
  \item
 \emph{orthogonal sum}:
\[
(A\perp B)_x=A_x\oplus B_x,\qquad x\in S_0\cup S_1.
\]

\end{itemize}

\noindent{\it Example 1.}
The problems of classifying, up to congruence, the representations over a skew field $K$ of the dographs
    \begin{equation}\label{2.8akg}
\xymatrix{
 {v}
 \ar@(ur,dr)@{-}^{\alpha}&
 }\\[2mm]
\end{equation}
\begin{equation}\label{2.8akd}
\xymatrix{
 {v} \ar@(ul,dl)@{-}_{\alpha}
 \ar@(ur,dr)@{-}^{\beta}&&
 {\alpha=\varepsilon \alpha^*,
 \quad \beta=\delta \beta^*,}
 }
\end{equation}
\begin{equation}\label{2.8akda}
\xymatrix{
 {v\,} \ar@(ul,dl)@{-}_{\alpha}
 \ar@<-0.4ex>@(u,r)@{-}^{\beta}
\ar@<-0.4ex>@(r,d)@{<->}^{\gamma}
 &&
  {\begin{matrix}
 \beta=\alpha^*\beta\alpha,\quad
 \beta=\varepsilon \beta^*\\
 \gamma\beta=1_v,\quad
 \beta\gamma=1_{v^*},
 \end{matrix}}
 }
\end{equation}
\begin{equation}\label{2.8akdag}
\xymatrix{
 {v\,} \ar@(ul,dl)@{-}_{\alpha}
 \ar@<-0.4ex>@(u,r)@{-}^{\beta}
\ar@<-0.4ex>@(r,d)@{<->}^{\gamma}
 &&
  {\begin{matrix}
 \beta\alpha=\alpha^*\beta,\quad
 \beta=\varepsilon \beta^*\\
 \gamma\beta=1_v,\quad
 \beta\gamma=1_{v^*}
 \end{matrix}}
 }
\end{equation}
(where $\varepsilon$ and $\delta$ are elements of the center of $K$, $\varepsilon \bar\varepsilon =\delta \bar\delta =1$)\footnote{The edge $\gamma$ and the
relations $\gamma\beta=1_v,$
$\beta\gamma=1_{v^*}$ ensure
the nonsingularity of forms
assigned to $\beta$.}, are the problems of classifying, respectively:
\begin{itemize}
  \item
sesquilinear forms over
$K$,

  \item
pairs of forms, the first
form is
$\varepsilon$-Hermitian and
the second is
$\delta$-Hermitian,

  \item
isometric operators on a
space with nondegenerative
$\varepsilon$-Hermitian form
(an operator $A$ is
\emph{isometric} for a form
$F(x,y)$ if
$F(Ax,Ay)=F(x,y)$), and
  \item
selfadjoint operators on a
space with nondegenerative
$\varepsilon$-Hermitian form
(an operator $A$ is
\emph{selfadjoint} for a form
$F(x,y)$ if
$F(x,Ay)=F(Ax,y)$).
\end{itemize}

\medskip

\noindent{\it Example 2.}
The problem of classifying the representations of a group $G$ by isometries of a nondegenerate $\varepsilon$-Hermitian form is presented by the dograph \eqref{2.8akda}, with the arrow $\alpha $ replaced by arrows $\alpha_1,\dots,\alpha _n$ (these being generators of $G$), and the relation $\beta = \alpha^*\beta \alpha$  replaced  by the relations $\beta_i = \alpha_i^*\beta_i \alpha_i$ ($1\le i\le n$) and the defining relations of $G$ (see \cite[Chap. 7, no. 2.6]{w.schar}).
\medskip

The rest of the paper has to do with the representations of the dographs \eqref{2.8akg}--\eqref{2.8akdag} over a field $K$ of characteristic not $2$ (these, as well as the representations of some other dographs, were announced in \cite{ser_disch,ser_dep}. Without loss of generality,
we assume that $\varepsilon,
\delta\in\{-1,1\}$ in the
case of the identity
involution on $K$,
and $\varepsilon= \delta=1$
in the case of nonidentity
involution (in the case of
nonidentity involution, an
$\varepsilon$-Hermitian form
can be made Hermitian by
multiplying by
$1+\bar{\varepsilon}$ if
$\varepsilon\ne -1$, and by
$a-\bar a\ne 0$ if
$\varepsilon\ne -1$).

For any polynomial
\[
f(x)=a_0x^n+a_1x^{n-1}+\dots
+a_n\in K[x]
\]
we define the polynomials
\begin{align*}\label{mau}
f^{\vee}(x)&:=\bar a_n^{-1}(\bar
a_nx^n+\dots+\bar a_1x+\bar
a_0)\quad\text{if } a_n\ne
0,\\
\bar f(x)&:=\bar a_0x^n+\bar
a_1x^{n-1}+\dots+\bar a_n.
\end{align*}
By the \emph{adjoint} of the matrix $A=[a_{ij}]$, we mean the matrix $A^*=[\bar a_{ji}]$) (this being the
matrix of the adjoint operator on the adjoint bases).

Every square matrix over
$K$ is similar to a
direct sum of {\it Frobenius
blocks}
\begin{equation*}\label{3.lfo}
\Phi=\begin{bmatrix} 0&&
0&-c_n\\1&\ddots&&\vdots
\\&\ddots&0&-c_2\\
0&&1& -c_1 \end{bmatrix},
\end{equation*}
whose characteristic
polynomials
\[
\chi_{\Phi}=x^n+c_1
x^{n-1}+\dots+ c_n
\]
are integer powers of
polynomials $p_{\Phi}(x)$
that are irreducible over
$K$.
For each Frobenius block
$\Phi$, denote by
$\sqrt[\displaystyle
*]{\Phi}$,\label{pagggs}
$\Phi_{\varepsilon}$, and
$\Phi_{(\varepsilon)}$
($\varepsilon=\pm 1$,
$\varepsilon=1$ for
nonidentity involution on
$K$) fixed
nonsingular matrices
satisfying, respectively, the
conditions\footnote{The matrix $\sqrt[\displaystyle
*]{\Phi}$ was denoted by $\widehat\Phi$ in [V. V. Sergeichuk,
Classification
problems  for systems
of forms and linear
mappings, {\it  Math.
USSR-Izv.} 31 (no. 3)
(1988) 481--501].}
\begin{align}\label{vrfmau1}
\sqrt[\displaystyle
*]{\Phi}&=(\sqrt[\displaystyle
*]{\Phi})^*\Phi,
   \\ \label{vrfmau2}
\Phi_{\varepsilon}
&=\Phi_{\varepsilon}^*,\quad
\Phi_{\varepsilon}\Phi=
\varepsilon
(\Phi_{\varepsilon}\Phi)^*,
  \\ \label{vrfmau3}
\Phi_{(\varepsilon)}
&=\varepsilon
\Phi_{(\varepsilon)}^*=
\Phi^*\Phi_{(\varepsilon)}\Phi.
\end{align}
Each of these matrices may do
not exist for some $\Phi$;
existence conditions and
explicit forms of these
matrices will be established
in Section 3.

The following lemma will be
employed in the construction
of the set
$\ind_0(\overline{{S}})$.

\begin{lemma}\label{lemaaay}
Let $S$ be a
dograph. If a
representation $\cal A$ of
the quiver $\overline S$ is isomorphic to a
selfadjoint representation, then
there exist a selfadjoint
representation $B$ and
an isomorphism $h: A\to B$ such that
$h_v=1$ for all vertices $v$
of $S$.
\end{lemma}

\begin{proof}
Let $f\colon {A}\to
{C}={C}^{\circ}$ be an isomorphism.
Define ${B}={B}^{\circ}$ and a congruence $g\colon {A}\to {B}$ as follows:
\[
g_u:=f_u^{-1},\quad
g_{u^*}:=f_u^*,\quad {B}_u:={A}_u,\quad {B}_{u^*}:={A}_u^*
\]
for each vertex $u$ of $S$, and
\[
{B}_{\alpha}:=g_v {C}_{\alpha}g_u^{-1}
\]
for each  arrow $\alpha\colon
u\to v$. Then $h:=gf\colon
{A}\to {B}$ is the
desired isomorphism.
\end{proof}

\subsection{Classification of
sesquilinear
forms}\label{s_ses}

\begin{lemma}\label{LEMMA 7}
Let $p(x) = p^{\vee}(x)$ be
an irreducible polynomial of
degree $2r$ or $2r+1$. Then
every stationary element of
the field
\begin{equation}\label{alft}
K(\kappa) = K[x]/p(x)K[x]
\end{equation}
with the involution
\begin{equation}\label{alfta}
f(\kappa)^{\circ} := \bar
f(\kappa^{-1})
\end{equation}
is uniquely representable in
the form $q(\kappa)$, where
\begin{equation}\label{ser13}
    q(x)=\bar a_rx^{-r}+\dots
    +a_0+\dots+a_rx^r
\end{equation}
$(a_0 = \bar a_0,a_1,\dots
a_r\in K)$, and when
$\deg(p(x)) = 2r$ the
following hold:
\begin{itemize}
  \item[\rm(a)]
 $a_r = 0$ if the involution
on $K$ is the
identity.
  \item[\rm(b)]
 $a_r=\bar a_r$ if the
involution on $K$ is
nonidentity and $p(0)\ne 1$.
  \item[\rm(c)]
 $a_r=-\bar a_r$ if the
involution on $K$ is
nonidentity and $p(0)=1$.
\end{itemize}
\end{lemma}

\begin{proof}
Suppose $\deg(p(x)) = 2r +1$.
The elements
$$\kappa^{-r},\dots,1,\dots,
\kappa^r$$ are linearly
independent over $K$.
Therefore all elements of the
form
$$a_{-r}\kappa^{-r}+\dots
    +a_0+\dots+a_r\kappa^r$$
are distinct. They are
stationary if and only if
$a_{-i} = \bar a_i$ for all
$i=0,1,\dots,r$.

Suppose $\deg(p(x)) = 2r$ and
the involution on $K$
is the identity. Then the
stationary elements of the
form
$$a_{r-1}\kappa^{-r+1}+\dots
    +a_0+\dots+a_{r-1}\kappa^{r-1}$$
are distinct and form a
vector space of dimension $r$
over $K$. But this is
the dimension over $K$ of the whole stationary
subfield of the field
$K(\kappa)$, and
therefore the subfield and
the vector space coincide.

Suppose $\deg(p(x)) = 2r$ and
the involution on $K$
is nonidentity. The equality
$p(x) = p^{\vee}(x)$ implies
that $\alpha\bar\alpha= 1$,
where $\alpha = p(0)$. Taking
any $a\ne\bar a\in K$
and putting
\[
\delta=
  \begin{cases}
    1 +
\bar\alpha & \text{if }
\alpha\ne 1,
 \\
    a-\bar a& \text{if }
    \alpha=1,
  \end{cases}
\]
we find that
$\delta\alpha=\bar\delta$.
The function $\pi(x): =
\delta x^{-r}p(x)$ has the
form
\[
 \pi(x)= c_{-r}x^{-r} +\dots+
 c_rx^r\qquad   (c_{-i}=\bar
 c_i).
\]
Using the equalities
$c_r=\delta$ and
$\delta\alpha=\bar\delta$, we
find that $c_r\ne \bar c_r$
if $\alpha\ne 1$, and $c_r\ne
-\bar c_r$ if $\alpha= 1$.

Let $q(x)$ be a function of
the form \eqref{ser13}. If
$q(\kappa) = 0$, then $q(x) =
a\pi(x)$, $a = \bar
a\in K$, and in view
of conditions (b) and (c) of
the lemma this is possible
only if $q(x) = 0$.
Consequently, the stationary
elements $q(\kappa)$ are
distinct and form a vector
space of dimension $2r$ over
the stationary subfield
$K_0$ of $K$.
But this is the dimension
over $K_0$ of the
whole stationary subfield of
$K(\kappa)$.
\end{proof}

Define the {\it skew sum} of
two matrices $A$ and $B$ as
follows:
\begin{equation*}\label{1.2a}
[A\diag B]=
\begin{bmatrix}0&B\\A &0
\end{bmatrix}.
\end{equation*}

\begin{theorem}\label{THEOREM
3}
\footnote{See also [R.A. Horn, V.V. Sergeichuk, Canonical matrices of bilinear and sesquilinear forms, Linear Algebra Appl. 428 (2008) 193--223; arXiv:0709.2408].}
Let $K$ be a
field of characteristic not
two with involution $($the
involution can be the
identity$)$. For  any
sesquilinear form on a
finite-dimensional vector
space over $ K$, there
exists a basis, in which the
matrix of the form is a
direct sum of matrices of the
three types:
\begin{itemize}
 \item [{\rm(i)}]
a singular Jordan block
$J_n(0)$;

 \item [{\rm(ii)}]
$A_{\Phi}=[\Phi\diag I]$,
where $\Phi$ is a nonsingular
Frobenius block for which
$\sqrt[\displaystyle
*]{\Phi}$ does not exist; the
block $\Phi$ and the identity
matrix $I$ have the same
size, and

 \item [{\rm(iii)}]
$\sqrt[\displaystyle
*]{\Phi}q(\Phi)$, where
$q(x)$ is a nonzero function
of the form \eqref{ser13}.
\end{itemize}
The summands are determined
to the following extent:
\begin{description}
  \item [Type (i)] uniquely.

  \item [Type (ii)]
up to replacement of the
block $\Phi$ by the block
$\Psi$ with $\chi_{\Psi}(x)
=\chi^{\vee}_{\Phi}(x)$.

  \item [Type (iii)]
up to replacement of the
whole group of summands
\[
\sqrt[\displaystyle
*]{\Phi}q_1(\Phi)
\oplus\dots\oplus
\sqrt[\displaystyle
*]{\Phi}q_s(\Phi)
\]
with the same $\Phi$ by
\[
\sqrt[\displaystyle
*]{\Phi}q'_1(\Phi)
\oplus\dots\oplus
\sqrt[\displaystyle
*]{\Phi}q'_s(\Phi)
\]
in which each $q'_i(x)$ is a
nonzero function of the form
\eqref{ser13} and the
Hermitian forms
\begin{gather}\label{777}
q_1(\kappa)x_1^{\circ}x_1+\dots+
q_s(\kappa)x_s^{\circ}x_s
\\\label{777s}
q'_1(\kappa)x_1^{\circ}x_1+\dots+
q'_s(\kappa)x_s^{\circ}x_s
\end{gather}
are equivalent over the
field\/ ${ K}[\kappa]={ K}[x]/p_{\Phi}{ K}[x]$
defined in \eqref{alft} with
the involution \eqref{alfta}.

In particular, if $ K$
is an algebraically closed
field with the identity
involution, then the summands
of type {\rm(iii)} can be
taken equal to
$\sqrt[\displaystyle
*]{\Phi}$. If $ K$ is
an algebraically closed field
with nonidentity involution,
or a real closed field,
then the summands of type
{\rm(iii)} can be taken equal
to $\pm\sqrt[\displaystyle
*]{\Phi}$. In these cases the
summands are then uniquely
determined by the
sesquilinear form.
\end{description}
\end{theorem}

\begin{proof}
We will study representations
of the dograph
\eqref{2.8akg}:
\[
\xymatrix{
 {v}
 \ar@(ur,dr)@{-}^{\alpha}&
 }\\[2mm]
\]

$1^{\circ}$ Let us describe
$\ind(\overline{{S}})$.
The dograph
${{S}}$ defines the quiver
\begin{equation}\label{ser14}
\overline{{S}}:\quad\xymatrix{
 {v}
  \ar@/^/@{->}[rr]^{\alpha}
 \ar@/_/@{->}[rr]_{\alpha^*} &&{v^*}
 }
\end{equation}
The representations of this
quiver, as well as morphisms
of the representations, will
be specified by pairs of
matrices. A representation is
a matrix pair
$(A_{\alpha},A_{\alpha^*})$
of the same size with entries
in $ K$. A morphism
$g\colon
(A_{\alpha},A_{\alpha^*})\to
(B_{\alpha},B_{\alpha^*})$ is
a matrix pair
$g=[G_v,G_{v^*}]$ (for
morphisms we use square
brackets) such that
\begin{equation}\label{ser15}
G_{v^*}A_{\alpha} =
B_{\alpha}G_v,\qquad
G_{v^*}A_{\alpha^*} =
B_{\alpha^*}G_v.
\end{equation}
The adjoint of a
representation is given by
\[
(A_{\alpha},A_{\alpha^*})^{\circ}:=
(B_{\alpha},B_{\alpha^*}),
\]
where
$B_{\alpha}:=A_{\alpha^*}^*$
and
$B_{\alpha^*}:=A_{\alpha}^*$

As shown by Kronecker (the
matrix pencil problem; see
\cite[Chap. XII]{gan}), the
set $\ind(\overline {{S}})$ consists of the
representations
\begin{equation}\label{ser16}
(N_1,N_2),\ (N_2^*,N_1^*),\
(\Phi, I_n),\, (I,J_n(0)),
\end{equation}
where $\Phi$ is an $n\times
n$ Frobenius block and
\begin{equation}\label{ser17}
N_1:=\begin{bmatrix}
1&0&&0\\&\ddots&\ddots&\\0&&1&0
\end{bmatrix},\quad
N_2:=\begin{bmatrix}
0&1&&0\\&\ddots&\ddots&\\0&&0&1
\end{bmatrix}.
\end{equation}

$2^{\circ}$. We describe
$\ind_0(\overline{{S}})$
and $\ind_1(\overline{{S}})$. By \eqref{ser15},
\[
(\Psi,I_n)\simeq
(\Phi,I_n)^{\circ}=(I_n,\Phi^*)
\]
if and only if $\Psi$ is
similar to $\Phi^{*-1}$,
i.e., if and only if
\begin{align*}
\chi_{\Psi}(x)&=
\det(xI_n-\Phi^{*-1})\\&=
\det(-\Phi^{*-1})\cdot
x^n\cdot\det(x^{-1}-\Phi^*)
=\chi_{\Phi}^{\vee}(x).
\end{align*}

Suppose the representation
$(\Phi, I_n)$ is isomorphic
to a selfadjoint
representation. By Lemma
\ref{lemaaay}, there exists
an isomorphism
\[
h=(I,H)\colon (\Phi,I)\to (A,
A^*).
\]
By \eqref{ser15}, $A = H\Phi$
and $A^* = H$. Then
$A=A^*\Phi$, and by
\eqref{vrfmau1} we can take
\begin{equation}\label{ser18}
h=(I,\sqrt[\displaystyle
*]{\Phi}^{\,*})\colon
(\Psi,I)\to
(\sqrt[\displaystyle
*]{\Phi}, \sqrt[\displaystyle
*]{\Phi}^{\,*}).
\end{equation}

Consequently, the set
$\ind_0(\overline{{S}})$
consists of the
representations $M_{\Phi} =
(\sqrt[\displaystyle
*]{\Phi},\sqrt[\displaystyle
*]{\Phi}^{\,*})$. The set
$\ind_1(\overline{{S}})$
consists of the
representations $(N_1, N_2)$
and $(\Phi,I)$, where $\Phi$
is a Frobenius block for
which $\sqrt[\displaystyle
*]{\Phi}$ does not exist; and
if $\Phi$ is nonsingular,
then it is determined up to
replacement by the Frobenius
block with characteristic
polynomial
$\chi^{\vee}_{\Phi}(x)$.

$3^{\circ}$. We describe the
orbits of the representations
in $\ind_0(\overline{{S}})$. Let $g= [G_1, G_2]\in
\End(M_{\Phi})$ and $h$ be
the isomorphism
\eqref{ser18}. Then
\[
h^{-1}gh =
[G_1,\sqrt[\displaystyle
*]{\Phi}^{\,*-1}G_2\sqrt[\displaystyle
*]{\Phi}^{\,*}]\colon
(\Phi,I)\to (\Phi,I);
\]
that is,
\[
G_1=\sqrt[\displaystyle
*]{\Phi}^{\,*-1}G_2\sqrt[\displaystyle
*]{\Phi}^{\,*},\qquad
 \Phi G_1=G_1\Phi.
\]
Since a matrix that commutes
with a Frobenius block is a
polynomial in this block, we
have
\begin{align*}
G_1 &=f(\Phi),\qquad
f(x)\in K[x],\\
G_2 &= \sqrt[\displaystyle
*]{\Phi}^{\,*}f(\Phi)
\sqrt[\displaystyle
*]{\Phi}^{\,*-1}=
f(\sqrt[\displaystyle
*]{\Phi}\,\Phi
\sqrt[\displaystyle
*]{\Phi}^{\,*-1})=
f({\Phi}^{\,*-1}).
\end{align*}

Consequently, the ring
$\End(M_{\Phi})$ consists of
matrix pairs
\[
g_f = [f(\Phi),
f(\Phi^{*-1})],\qquad
f(x)\in K[x],
\]
with involution
\[
g_f^{\circ} = [\bar
f(\Phi^{-1}), \bar
f(\Phi^{*})].
\]
By Lemma \ref{lemser1a}, its
radical $R$ consists of the
pairs $g_f$ for which
$f(x)\in p_{\Phi(x)} K[x]$. Hence the field
$T(M_{\Phi}) =
\End(M_{\Phi})/R$ can be
identified with the field
$ K[\kappa] =  K[x]/p_{\Phi}(x)$ with
involution $f(\kappa)^{\circ}
= \bar f(\kappa^{-1})$.

Under this identification, a
stationary element
$q(\kappa)\ne 0$ of the field
$ K[\kappa]$ (where
$q(x)$ is a function of form
\eqref{ser13}) corresponds to
the coset in the quotient
ring $\End(M_{\Phi})/R$ that
contains the selfadjoint
automorphism $[q(\Phi),
q(\Phi^{*-1})]$. By
\eqref{ser02a}, the
representations
\[
M_{\Phi}^{q(\kappa)}=(\sqrt[\displaystyle
*]{\Phi}q(\Phi),\sqrt[\displaystyle
*]{\Phi}^{\,*}q(\Phi))
\]
constitute the orbit of the
representation $M_{\Phi}$.

$4^{\circ}$. We now apply
Theorems \ref{THEOREM 1} and \ref{THEOREM 2}. Each
selfadjoint representation
$(A, A^*)$ of the quiver
(\ref{ser14}) corresponds, in
a one-to-one manner, to the
representation of the
dograph
\eqref{2.8akg} given by the
matrix $A$. In particular,
the representation $(A, B)^+$
of the quiver (\ref{ser14})
corresponds to the
representation $[A\diag B^*]$
(see \eqref{1.2a}) of the
dograph. From
Theorem \ref{THEOREM 1} and
items $2^{\circ}$ and
$3^{\circ}$ above, it follows
that every representation of
the dograph
\eqref{2.8akg} is congruent
to an orthogonal sum of representations
of the form $[N_1\diag
N_2^*]$, $[\Phi \diag I]$
if $\sqrt[\displaystyle
*]{\Phi}$ does not exist,
and $\sqrt[\displaystyle
*]{\Phi}f(\Phi)$.

 Let us prove that the
representations $[N_1\diag
N_2^*]$  and $[J_n(0) \diag
I_n]$ are congruent to a
singular Jordan block. We
show that each matrix
$[N_1\diag N_2^T]$ or
$[J_m(0)\diag I_m]$ can be
obtained by simultaneous
permutations of rows and
columns of a singular Jordan
block. The units of $J_n(0)$
are disposed at the places $
(1,2),\
(2,3),\,\dots,\,(n-1,n); $ it
suffices to prove that there
is a permutation $f$ on
$\{1,2,\dots,n\}$ such that
\[
(f(1),f(2)),\ \ (f(2),f(3)),\
\dots,\ (f(n-1),f(n))
\]
are the positions of the unit
entries in $[N_1\diag N_2^T]$
if $n=2m-1$ or in
$[J_m(0)\diag I_m]$ if
$n=2m$. This becomes clear if
we arrange the positions of
the unit entries in the
$(2m-1)\times(2m-1)$ matrix
$$
[N_1\diag N_2^T]=\left[
\begin{tabular}{cccc|ccc}
 &&& & 0&&0\\
 &0&& & 1&$\ddots$&\\
 &&& & &$\ddots$&0\\
 &&& & 0&&1
 \\ \hline
 1&0&&0 &&&\\
 &$\ddots$&$\ddots$& &&0&\\
 0&&1&0 &&&
\end{tabular}\right]
$$
as follows:
\begin{multline*}
(m,2m-1),\, (2m-1,m-1),\,
(m-1,2m-2),\\
(2m-2,m-2),\dots,(2,m+1),\,
(m+1,1),
\end{multline*}
and the positions of the unit
entries in the $2m\times 2m$
matrix $[J_m(0)\diag I_m]$ as
follows:
\begin{multline*}
(1,m+1),\,(m+1,2),\,(2,m+2),\\
(m+2,3),
\dots,(2m-1,m),\,(m,2m).
\end{multline*}

This proves the first
assertion of Theorem
\ref{THEOREM 3} (concerning
existence of a basis). The
remaining assertions follow
from Theorems \ref{THEOREM 1} and \ref{THEOREM 2}.
\end{proof}

\noindent{\it Remark.}%
\footnote{This statement was proved in
[V.V. Sergeichuk, The
canonical form of the
matrix of a bilinear
form over an
algebraically closed
field of
characteristic 2, {\it
Math. Notes} 41 (1987)
441--445.]}
 It can be shown that over an algebraically closed field of characteristic $2$, there exists for any bilinear form a basis in which its matrix is a direct sum
\[
[\Phi_1\diag I]\oplus\dots\oplus [\Phi_r\diag I]\oplus\sqrt[\displaystyle
*]{\Psi}_1\oplus\dots\oplus  \sqrt[\displaystyle
*]{\Psi}_t\oplus J_{n_1}(0)\oplus\dots\oplus J_{n_s}(0),
\]
where the $\Phi_i$ and $\Psi_j$ are nonsingular Jordan blocks and $\Phi_i\ne \Psi_j$ for all $i,j$.
This direct sum is uniquely determined by the bilinear form up to permutation of the summands and replacement of the eigenvalue $\lambda $ in a block $\Phi_i$ by $\lambda ^{-1}$. The matrix $\sqrt[\displaystyle
*]{\Psi}$ exists if and only if the matrix ${\Psi}$ is of odd size with eigenvalue $1$.

\subsection{Classification of
pairs of Hermitian
forms}\label{s_her}

\begin{lemma}\label{LEMMA 8}
Let K be a field with the identity
involution, and suppose $A =
\varepsilon A$ and $A\Phi =
\delta(A\Phi)^*$ for a
nonsingular matrix $A$ and a
Frobenius block $\Phi$. Then
either $\varepsilon = 1$ or
$\delta = 1$. If
$\chi_{\Phi}= x^n$, then
$\varepsilon = 1$ for $n$ odd
and $\delta = 1$ for $n$
even.
\end{lemma}

\begin{proof}
Let $A = [a_{ij}]$ be
$n$-by-$n$. Since
multiplication by a Frobenius
block moves the columns of
this matrix to the left, we
have $A\Phi = [a_{i,j+1}]$
(the entries $a_{i,n+1}$ are
defined by this equality).
The relations $$A=\varepsilon
A^*,\qquad A\Phi =
\delta(A\Phi)^*$$ can then be
written
\begin{equation}\label{ser19}
a_{ij} = \varepsilon
a_{ji},\qquad
a_{i,j+1}=\delta a_{j,i+1}.
\end{equation}
Consequently,
\[
a_{ij} = \varepsilon\delta
a_{i-1,j+1}=
(\varepsilon\delta)^{-i}b_{i+j},\qquad
b_2,\dots,b_{2n}\in K.
\]
Putting $i= j$ in
\eqref{ser19}, we find that
$b_{2i} = 0$ if
$\varepsilon\ne ~ 1$, and
$b_{2i+1}= 0$ if $\delta\ne
1$. Since $A\ne 0$, this
implies either $\varepsilon =
1$ or $\delta= 1$. If
$\chi_{\Phi}(x)= x^n$, then
the formula $A\Phi
=[b_{i+j+1}]$ implies
\[
b_{n+2}=b_{n+3}=\dots=0;
\]
and since $A=[b_{i+j}]$ is
nonsingular, this means
$b_{n+1}\ne 0$, and therefore
$\varepsilon= 1$ for $n$ odd,
$\delta = 1$ for $n$ even.
\end{proof}

For any matrices $A, B, C, D$
we define
\[
(A,B)\oplus (C,D)=(A\oplus
C,B\oplus D),\qquad (A,B)C =
(AC,BC).
\]

\begin{theorem}\label{Theorem 4}
Let $F_1$ and $F_2$ be
$\varepsilon$- and
$\delta$-Hermitian forms,
respectively, in a finite
dimensional vector space over
a field $ K$ of
characteristic $\ne 2$
$(\varepsilon = \pm 1$,
$\delta = \pm 1$,
$\varepsilon \ge \delta$, and
$\varepsilon = \delta=1$ for
nonidentity involution on
$ K)$. Then there
exists a basis in which the
pair $(F_1,F_2)$ is given by
a direct sum of matrix pairs
of the following types:
\begin{itemize}
  \item[\rm(i)]
$([N_1\diag \varepsilon
N_1^*],\, [N_2\diag \delta
N_2^*])$, where $N_1$ and
$N_2$ are defined in
\eqref{ser17}.

  \item[\rm(ii)]
$([I_n\diag \varepsilon I_n],\,
[\Phi\diag \delta \Phi^*])$,
where $\Phi$ is an $n\times n$ Frobenius
block such that
$\Phi_{\delta}$ $($see
\eqref{vrfmau2}$)$ does not
exist if $\varepsilon=1$.

  \item[\rm(iii)]
$A_{\Phi}^{f(x)}:=
(\Phi_{\delta},
\Phi_{\delta}\Phi)f(\Phi)$,
where $\varepsilon=1$, $0\ne
f(x)=\bar f(\delta
x)\in K[x]$, and
$\deg(f(x))<\deg(p_{\Phi}(x))$.

  \item[\rm(iv)]
$([J_n(0)\diag \varepsilon
J_n(0)^*], [I_n\diag
(-I_n)])$, where $\delta=-1$,
and $n$ is odd if
$\varepsilon=1$.

  \item[\rm(v)] $ B_n^a:=$
\begin{equation}\label{ser20}
\!\!\!\!\!\!\!\!\!\!\!
\left(\!
 \begin{bmatrix}
0&&&&&1&0\\
&&&&\delta&\cdot&\\
&&&1&\cdot&&\\
&&\delta&\cdot&&&\\
&\cdot&\cdot&&&&\\
\cdot&\cdot&&&&&\\
0&&&&&&0
\end{bmatrix}\!a,
\begin{bmatrix}
0&&&&&&1\\
&&&&&\delta&\\
&&&&1&&\\
&&&\delta&&&\\
&&\cdot&&&&\\
&\cdot&&&&&\\
&&&&&&0
\end{bmatrix}\!a
 \!\right)
\end{equation}
where the matrices are
$n$-by-$n$, $\varepsilon=1$,
$0\ne a=\bar a\in K$,
and $n$ is even if
$\delta=-1$.

The summands are determined
to the following extent:
\begin{description}
  \item [Type (i)]
uniquely.

  \item [Type (ii)]
up to replacement of $\Phi$
by $\Psi$ with
$\chi_{\Psi}(x)
=\pm\bar\chi_{\Phi}(\varepsilon
\delta x)$.

  \item [Type (iii)]
up to replacement of the
whole group of summands
\[
A_{\Phi}^{f_1(x)}
\oplus\dots\oplus
A_{\Phi}^{f_s(x)}
\]
with the same $\Phi$ by
\[
A_{\Phi}^{g_1(x)}
\oplus\dots\oplus
  A_{\Phi}^{g_s(x)}
\]
such that the Hermitian forms
\begin{gather*}\label{777aki}
f_1(\omega)x_1^{\circ}x_1+\dots+
f_s(\omega)x_s^{\circ}x_s,
\\\label{777sopja}
g_1(\omega)x_1^{\circ}x_1+\dots+
g_s(\omega)x_s^{\circ}x_s
\end{gather*}
are equivalent over the
field\/ ${ K}[\omega]={ K}[x]/p_{\Phi}{ K}[x]$
with involution
$f(\omega)^{\circ}= \bar
f(\delta\omega)$.

  \item [Type (iv)]
uniquely.

  \item [Type (v)]
up to replacement of the
whole group of summands
\[
B_n^{a_1} \oplus\dots\oplus
B_n^{a_s}
\]
with the same $n$ by
\[
B_n^{b_1} \oplus\dots\oplus
B_n^{b_s}
\]
such that the Hermitian forms
\begin{gather*}\label{777a}
a_1x_1^{\circ}x_1+\dots+
a_sx_s^{\circ}x_s,
\\\label{777sa}
b_1x_1^{\circ}x_1+\dots+
b_sx_s^{\circ}x_s
\end{gather*}
are equivalent over\/
${ K}$.
\end{description}
\end{itemize}
\end{theorem}

\begin{proof}
We will study representations
of the dograph
\eqref{2.8akd}:
\begin{equation*}\label{2.8akdg}
{{S}}:\quad\xymatrix{
 {v} \ar@(ul,dl)@{-}_{\alpha}
 \ar@(ur,dr)@{-}^{\beta}&&
 {\alpha=\varepsilon \alpha^*,
 \quad \beta=\delta \beta^*.}
 }
\end{equation*}

$1^{\circ}$ Let us describe
$\ind(\overline{{S}})$.
The dograph
${{S}}$ defines the quiver
\begin{equation*}
\overline{{S}}:\quad\xymatrix{
 {v}
 \ar@/^1.5pc/@{->}[rr]^{\alpha}
  \ar@/^/@{->}[rr]^{\alpha^*}
 \ar@/_/@{->}[rr]_{\beta}
 \ar@/_1.5pc/@{->}[rr]_{\beta^*} &&{v^*}
 },\qquad
 \alpha=\varepsilon\alpha^*,\
 \ \beta=\delta\beta^*.
\end{equation*}
The representations of this
quiver will be specified by
pairs of matrices
$(A_{\alpha}, A_{\beta})$ of
the same size; then
$A_{\alpha^*}= \varepsilon
A_{\alpha}$ and $A_{\beta^*}=
\delta A_{\beta}$. The
adjoint representation is
given by
\[
(A_{\alpha},
A_{\beta})^{\circ} =
(\varepsilon A_{\alpha}^*,
\delta A_{\beta}^*).
\]
The set $\ind(\overline{{S}})$ consists of the
representations
\[
(N_1,N_2),\ (N_1^*,N_2^*),\
(I_n,\Phi),\, (J_n(0),I_n)
\]
(which we prefer now to the
set \eqref{ser16}).

$2^{\circ}$. We describe
$\ind_0(\overline{{S}})$
and $\ind_1(\overline{{S}})$. It is obvious that
\[
(I,\Psi)\simeq
(I,\Phi)^{\circ}=
(\varepsilon I, \delta
\Phi^*)
\]
if and only if $\Psi$ is
similar to $\varepsilon
\delta \Phi^*$, i.e., if and
only if $$\chi_{\Psi}(x)
=\pm\bar\chi_{\Phi}(\varepsilon
\delta x).$$

Suppose $(I,\Phi)$ is
isomorphic to a selfadjoint
representation. By Lemma
\ref{lemaaay}, there exists
an isomorphism
\[
 h=[I, H]\colon
 (I,\Phi)\to (A, B) =
  (A, B)^{\circ}.
\]
Then
\[
A = H,\quad B = H\Phi,\quad A
= \varepsilon A^*, \quad B=
\delta B^*;
\]
i.e., $$A = \varepsilon
A^*,\qquad A\Phi = \delta
(A\Phi)^*.$$ Since
$\varepsilon \ge\delta$, we
have by Lemma \ref{LEMMA 8}
that $\varepsilon = 1$, and
by \eqref{vrfmau2},
\begin{equation}\label{ser21}
 h=[I, \Phi_{\delta}]\colon
 (I,\Phi)\to (\Phi_{\delta},
 \Phi_{\delta}\Phi).
\end{equation}

Similarly, if
\[
  (J_n(0),I_n)\simeq (A, B) =
  (A, B)^{\circ},
\]
then \[B = \delta B^*,\qquad
BJ_n(0) = \varepsilon
(BJ_n(0))^*;\] by Lemma
\ref{LEMMA 8}, $\varepsilon =
1$, and $n$ is even if
$\delta = -1$. It is easily
verified that
$$(J_n(0),I_n)\simeq B_n,$$
where $B_n = B_n^1$ is of the
form \eqref{ser20}.

Consequently, the set
$\ind_0(\overline{{S}})$
is empty if $\varepsilon =
-1$, and consists of the
representations $$A_{\Phi} =
(\Phi_{\delta},
\Phi_{\delta}\Phi)$$ and
$B_n$ (where $n$ is even when
$\delta= -1$) if $\varepsilon
= 1$.

The set
$\ind_1(\overline{{S}})$
consists of the following
representations:
\begin{itemize}
  \item
$(N_1,N_2)$

  \item
$(I,\Phi)$, where
$\Phi_{\delta}$ does not
exist if $\varepsilon = 1$,
and $\chi_{\Psi}(x)$ is
determined up to replacement
by
$\bar\chi_{\Phi}(\varepsilon
\delta x)$.

  \item
$(J_n(0), I_n)$, where
$\delta = -1$, and $n$ is odd
if $\varepsilon= 1$.
\end{itemize}

$3^{\circ}$.  We describe the
orbits of the representations
in $\ind_0(\overline{{S}})$. Let $$g = [G_1, G_2]\in
\End(A_{\Phi}),$$ and $h$ be
the isomorphism
\eqref{ser21}. Then
\[
h^{-1}gh = [G_1,
\Phi_{\delta}^{-1}G_2
\Phi_{\delta}]\colon (I,
\Phi)\to (I, \Phi);
\]
i.e., $$G_1 =
\Phi_{\delta}^{-1}G_2
\Phi_{\delta},\qquad\Phi G_1
= G_1\Phi.$$ Since $G_1$
commutes with $\Phi$, we have
$$G_1= f(\Phi),\qquad f(x)\in
 K[x],$$ and by
\eqref{vrfmau2},
\[
G_2 = \Phi_{\delta}f(\Phi)
\Phi_{\delta}^{-1} =
f(\Phi_{\delta}\Phi
\Phi_{\delta}^{-1})=
f(\delta\Phi^*).
\]

Consequently, the ring
$\End(A_{\Phi})$ consists of
the matrix pairs
\[
g_f = [f(\Phi),
f(\delta\Phi^*)],\qquad
f(x)\in K[x],
\]
with involution $$g_f^{\circ}
= [\bar
f(\delta\Phi),f(\Phi)^*].$$
Hence the field $$T(A_{\Phi})
= \End(A_{\Phi})/R$$ can be
identified with the field
$$ K[\omega] =  K[x]/p_{\Phi}(x) K[x],$$ with involution
$f(\omega)^{\circ} = \bar
f(\delta \omega)$. The set of
representations
$$A_{\Phi}^{f(\omega)} =
A_{\Phi}f(\Phi),$$ where
$$0\ne f(x)=\bar f(\delta
x)\in K[x],\qquad
\deg(f(x))<\deg(p_{\Phi}(x)),$$
is the orbit of the
representation $A_{\Phi}$.

Similarly, $T(B_n)$ can be
identified with the field
$ K$, and the set of
representations of the form
$B_na$, where $0\ne a=\bar
a\in K$, is the orbit
of the representation $B_n$.

$4^{\circ}$. From
$2^{\circ}$,  $3^{\circ}$,
and Theorem \ref{THEOREM 1},
the proof of Theorem
\ref{Theorem 4} now follows.
\end{proof}

\subsection{Classification of
isometric
operators}\label{s_iso}

\begin{theorem}\label{Theorem 5}
\footnote{See also [V.V. Sergeichuk, Canonical matrices of isometric operators on indefinite inner product spaces, Linear Algebra Appl. 428 (2008) 154--192; arXiv:0710.0933].}
Let $A$ be an isometric
operator on a
finite-dimensional vector
space with nondegenerate
$\varepsilon$-Hermitian form
$F$ over a field $ K$
of characteristic not $2$. Then
there exists a basis in which
the pair $(A, F)$ is given by
a direct sum of matrix pairs
of the following types:
\begin{itemize}
  \item[\rm(i)]
$(\Phi\oplus\Phi^{*-1},
[I_n\diag \varepsilon I_n])$, where
$\Phi$ is a nonsingular $n\times n$
Frobenius block for which
$\Phi_{(\varepsilon)}$ $($see
\eqref{vrfmau3}$)$ does not
exist.

  \item[\rm(ii)]
$A^{q(x)}_{\Phi}=(\Phi,
\Phi_{(\varepsilon)}q(\Phi))$,
where $q(x)\ne 0$ is of the
form \eqref{ser13}.
\end{itemize}
The summands are determined
to the following extent:
\begin{description}
  \item [Type (i)]
up to replacement of $\Phi$
by $\Psi$ with
$\chi_{\Psi}(x)
=\chi_{\Phi}^{\vee}(x)$.

  \item [Type (ii)]
up to replacement of the
whole group of summands
\[
A_{\Phi}^{q_1(x)}
\oplus\dots\oplus
A_{\Phi}^{q_s(x)}
\]
with the same $\Phi$ by
\[
A_{\Phi}^{q'_1(x)}
\oplus\dots\oplus
  A_{\Phi}^{q'_s(x)}
\]
such that the Hermitian forms
\begin{gather*}
q_1(\omega)x_1^{\circ}x_1+\dots+
q_s(\omega)x_s^{\circ}x_s,
\\
q'_1(\omega)x_1^{\circ}x_1+\dots+
q'_s(\omega)x_s^{\circ}x_s
\end{gather*}
are equivalent over the
field\/ ${ K}[\kappa]={ K}[x]/p_{\Phi}{ K}[x]$
with involution
$f(\kappa)^{\circ}= \bar
f(\kappa^{-1})$.
\end{description}
\end{theorem}

\begin{proof}
We will study representations
of the dograph
\eqref{2.8akda}:
\begin{equation*}\label{2.8dakda}
{{S}}:\quad\xymatrix{
 {v\,} \ar@(ul,dl)@{-}_{\alpha}
 \ar@<-0.4ex>@(u,r)@{-}^{\beta}
\ar@<-0.4ex>@(r,d)@{<->}^{\gamma}
 &&
  {\begin{matrix}
 \beta=\alpha^*\beta\alpha,\quad
 \beta=\varepsilon \beta^*,\\
 \gamma\beta=1_v,\quad
 \beta\gamma=1_{v^*}.
 \end{matrix}}
 }
\end{equation*}

$1^{\circ}$ Let us describe
$\ind(\overline{{S}})$.
The dograph
${{S}}$ defines the quiver
$\overline{{S}}:$
\begin{equation*}
\xymatrix{
 {v}\ar@(ul,dl)@{->}_{\alpha}
 \ar@/^1.5pc/@{->}[rr]^{\beta}
  \ar@/^/@{->}[rr]^{\beta^*}
 \ar@/^/@{->}[rr];[]^{\gamma}
 \ar@/^1.5pc/@{->}[rr];[]^{\gamma^*} &&{v^*}
\ar@(ur,dr)@{->}^{\alpha^*}
 }\qquad
 {\begin{matrix}
 \beta=\alpha^*\beta\alpha,\quad
 \beta=\varepsilon \beta^*,\\
 \gamma\beta=1_v,\quad
 \beta\gamma=1_{v^*},\\
 \gamma^*\beta^*=1_v,\quad
 \beta^*\gamma^*=1_{v^*}.
 \end{matrix}}
\end{equation*}

The representations of this
quiver will be specified by
triples of square matrices
$(A_{\alpha}, A_{\beta},
A_{\alpha^*})$ of the same
size, where $A_{\beta}$ is
nonsingular and $$A_{\beta} =
A_{\alpha^*}A_{\beta}A_{\alpha},$$
and then
\[
A_{\beta^*}=\varepsilon^{-1}
A_{\beta},\quad A_{\gamma}=
A_{\beta}^{-1},\quad
A_{\gamma^*}=\varepsilon
A_{\beta}^{-1}.
\]
The adjoint representation is
given by $$(A,B,C)^{\circ} =
(C^*,\varepsilon B^*,A^*).$$

Every representation of the
quiver is isomorphic to one
of the form $(A, I, A^{-1})$.
The set $\ind(\overline{{S}})$ consists of the
representations $(\Phi, I,
\Phi^{-1})$, where $\Phi$ is
a Frobenius block.

$2^{\circ}$. We describe
$\ind_0(\overline{{S}})$
and $\ind_1(\overline{{S}})$. It is obvious that
\[
(\Psi,I,\Psi^{-1})\simeq
(\Phi,I,\Phi^{-1})^{\circ} =
(\Phi^{*-1},\varepsilon
I,\Phi^{*})
\]
if and only if $\Psi$ is
similar to $\Phi^{*-1}$,
i.e., if and only if
$\chi_{\Psi}(x)
=\chi_{\Phi}^{\vee}(x)$.

Suppose $(\Phi,I,\Phi^{-1})$
is isomorphic to a
selfadjoint representation.
By Lemma \ref{lemaaay}, there
exists an isomorphism
\[
 h=[I, H]\colon
 (\Phi,I,\Phi^{-1})\to (A,B,A^*),
 \qquad B=\varepsilon B^*.
\]
Then
\[
A = \Phi,\quad B = H,\quad
A^*H = H\Phi^{-1},\quad B
=\varepsilon B^*;
\]
i.e.,
\[
A = \Phi,\quad B =\varepsilon
B^*=\Phi^*B\Phi.
\]
By (\ref{vrfmau3}),
\[
h=[I,\Phi_{(\varepsilon)}]\colon
(\Phi,I,\Phi^{-1})\to
(\Phi,\Phi_{(\varepsilon)},\Phi^*).
\]

Consequently, the set
$\ind_0(\overline{{S}})$
consists of the
representations $$A_{\Phi} =
(\Phi,\Phi_{(\varepsilon)},\Phi^*).$$
The set
$\ind_1(\overline{{S}})$
consists of the
representations
$(\Phi,I,\Phi^{-1})$, in
which $\Phi$ is a Frobenius
block such that
$\Phi_{(\varepsilon)}$ does
not exist and
$\chi_{\Phi}(x)$ is
determined up to replacement
by $\chi_{\Phi}^{\vee}(x)$.

$3^{\circ}$. We describe the
orbits of the representations
in $\ind_0(\overline{{S}})$. Let $$g=[G_1, G_2]\in
\End(A_{\Phi}).$$ Then
\[
\Phi G_1 = G_1\Phi,\quad
\Phi_{(\varepsilon)} G_1 =
G_2\Phi_{(\varepsilon)},\quad
\Phi^* G_2 = G_2\Phi^*.
\]
Since $G_1$ commutes with the
Frobenius block, we have
\[
G_1 = f(\Phi)\ \ (f(x)\in
 K[x]),\quad G_2 =
\Phi_{(\varepsilon)} f(\Phi)
\Phi_{(\varepsilon)}^{-1}=
f(\Phi^{*-1}).
\]

Consequently, the algebra
$\End(A_{\Phi})$ consists of
the matrix pairs
\[
[f(\Phi),f(\Phi^{*-1})],\qquad
f(x)\in K[x],
\]
with involution
\[
[f(\Phi),f(\Phi^{*-1})]^{\circ}=
[\bar
f(\Phi^{-1}),f(\Phi)^*].
\]
The field $T(A_{\Phi})$ can
be identified with the field
$${ K}[\kappa]={ K}[x]/p_{\Phi}{ K}[x]$$ with involution
$f(\kappa)^{\circ}= \bar
f(\kappa^{-1})$.

Let $q(\kappa)$ (where
$q(x)\ne 0$ is of the form
(\ref{ser13})) be a
stationary element of this
field. The representations
$$A_{\Phi}^{q(\kappa)} =
(\Phi,\Phi_{(\varepsilon)}q(\Phi))$$
constitute the orbit of the
representation $A_{\Phi}$.

$4^{\circ}$. From
$2^{\circ}$,  $3^{\circ}$,
and Theorem \ref{THEOREM 1},
the proof of Theorem
\ref{Theorem 5} now follows.
\end{proof}

\subsection{Classification of
selfadjoint operators}
\label{s_self}

\begin{theorem}\label{Theorem 6}
Let $A$ be a selfadjoint
operator on a
finite-dimensional vector
space with nondegenerate
$\varepsilon$-Hermitian form
$F$ over a field $ K$
of characteristic not $2$
$(\varepsilon = \pm 1$;
$\varepsilon= 1$ for
nonidentity involution on
$ K)$. Then there
exists a basis in which the
pair $(A, F)$ is given by a
direct sum of matrix pairs of
the following types:
\begin{itemize}
  \item[\rm(i)]
$(\Phi\oplus\Phi^{*}, [I_n\diag
\varepsilon I_n])$, where $\Phi$
is an $n\times n$ Frobenius block and if
$\varepsilon =1$ then
$\Phi_1$ $($see
\eqref{vrfmau2}$)$ does not
exist.

  \item[\rm(ii)]
$A^{f(x)}_{\Phi}=(\Phi,
\Phi_1f(\Phi))$, where
$\varepsilon =1$, $0\ne
f(x)=\bar f(x)\in K[x]$, and
$\deg(f(x))<\deg(p_{\Phi}(x))$.
\end{itemize}
The summands are determined
to the following extent:
\begin{description}
  \item [Type (i)]
up to replacement of $\Phi$
by $\Psi$ with
$\chi_{\Psi}(x)
=\bar\chi_{\Phi}(x)$.

  \item [Type (ii)]
up to replacement of the
whole group of summands
\[
A_{\Phi}^{f_1(x)}
\oplus\dots\oplus
A_{\Phi}^{f_s(x)}
\]
with the same $\Phi$ by
\[
A_{\Phi}^{g_1(x)}
\oplus\dots\oplus
  A_{\Phi}^{g_s(x)}
\]
such that the Hermitian forms
\begin{gather*}
f_1(\omega)x_1^{\circ}x_1+\dots+
f_s(\omega)x_s^{\circ}x_s,
\\
g_1(\omega)x_1^{\circ}x_1+\dots+
g_s(\omega)x_s^{\circ}x_s
\end{gather*}
are equivalent over the
field\/ ${ K}[\omega]={ K}[x]/p_{\Phi}{ K}[x]$
with involution
$f(\omega)^{\circ}= \bar
f(\omega)$.
\end{description}
\end{theorem}

\begin{proof}
We will study representations
of the dograph
\eqref{2.8akdag}:
\[
\xymatrix{
 {v\,} \ar@(ul,dl)@{-}_{\alpha}
 \ar@<-0.4ex>@(u,r)@{-}^{\beta}
\ar@<-0.4ex>@(r,d)@{<->}^{\gamma}
 &&
  {\begin{matrix}
 \beta\alpha=\alpha^*\beta,\quad
 \beta=\varepsilon \beta^*,\\
 \gamma\beta=1_v,\quad
 \beta\gamma=1_{v^*}.
 \end{matrix}}
 }
\]

$1^{\circ}$ Let us describe
$\ind(\overline{{S}})$.
The dograph
${{S}}$ defines the quiver
$\overline{{S}}:$
\begin{equation*}
\xymatrix{
 {v}\ar@(ul,dl)@{->}_{\alpha}
 \ar@/^1.5pc/@{->}[rr]^{\beta}
  \ar@/^/@{->}[rr]^{\beta^*}
 \ar@/^/@{->}[rr];[]^{\gamma}
 \ar@/^1.5pc/@{->}[rr];[]^{\gamma^*} &&{v^*}
\ar@(ur,dr)@{->}^{\alpha^*}
 }\qquad
 {\begin{matrix}
\beta\alpha=\alpha^*\beta,\quad
 \beta=\varepsilon \beta^*,\\
  \gamma\beta=1_v,\quad
 \beta\gamma=1_{v^*},\\
 \gamma^*\beta^*=1_v,\quad
 \beta^*\gamma^*=1_{v^*}.
 \end{matrix}}
\end{equation*}

The representations of this
quiver will be specified by
triples of square matrices
$(A_{\alpha}, A_{\beta},
A_{\alpha^*})$ of the same
size, where $A_{\beta}$ is
nonsingular and
$$A_{\beta}A_{\alpha} =
A_{\alpha^*}A_{\beta^*}.$$
 The
adjoint representation is
given by $$(A,B,C)^{\circ} =
(C^*,\varepsilon B^*,A^*).$$

Every representation of the
quiver is isomorphic to one
of the form $(A, I, A)$. The
set $\ind(\overline{{S}})$ consists of the
representations $(\Phi, I,
\Phi)$, where $\Phi$ is a
Frobenius block.

$2^{\circ}$. We describe
$\ind_0(\overline{{S}})$
and $\ind_1(\overline{{S}})$. It is obvious that
\[
(\Psi,I,\Psi)\simeq
(\Phi,I,\Phi)^{\circ} =
(\Phi^{*},\varepsilon
I,\Phi^{*})
\]
if and only if $\Psi$ is
similar to $\Phi^{*}$, i.e.,
if and only if
$\chi_{\Psi}(x)
=\bar\chi_{\Phi}(x)$.

Suppose $(\Phi,I,\Phi)$ is
isomorphic to a selfadjoint
representation. By Lemma
\ref{lemaaay}, there exists
an isomorphism
\[
 h=[I, H]\colon
 (\Phi,I,\Phi)\to (A,B,A^*),
 \qquad B=\varepsilon B^*.
\]
Then
\[
A = \Phi,\quad B = H,\quad
A^*H = H\Phi,\quad B
=\varepsilon B^*;
\]
i.e.,
\[
B =\varepsilon B^*,\quad
B\Phi =\Phi^* B=\varepsilon
(B\Phi)^*.
\]
By Lemma \ref{LEMMA 8},
$\varepsilon =1$ and we can
take $B=\Phi_1$ (see
\eqref{vrfmau2}).

Consequently, the set
$\ind_0(\overline{{S}})$
is empty if $\varepsilon =
-1$, and consists of the
representations $$A_{\Phi} =
(\Phi,\Phi_1,\Phi)$$ if
$\varepsilon = 1$. The set
$\ind_1(\overline{{S}})$
consists of the
representations
$(\Phi,I,\Phi)$, in which
$\Phi$ is a Frobenius block
such that $\chi_{\Phi}(x)$ is
determined up to replacement
by $\bar\chi_{\Phi}(x)$ and
if $\varepsilon = 1$ then
$\Phi_1$ does not exist.

$3^{\circ}$. We describe the
orbits of the representations
in $\ind_0(\overline{{S}})$. Let $$g=[G_1, G_2]\in
\End(A_{\Phi}).$$ Then
\[
\Phi G_1 = G_1\Phi,\quad
\Phi_1 G_1 = G_2\Phi_1,\quad
\Phi^* G_2 = G_2\Phi^*.
\]
Since $G_1$ commutes with the
Frobenius block, we have
\begin{align*}
G_1 &= f(\Phi)\qquad (f(x)\in
 K[x]),\\ G_2 &=
\Phi_1 f(\Phi) \Phi_1^{-1}=
f(\Phi_1 \Phi \Phi_1^{-1})=
f(\Phi^{*}).
\end{align*}

Consequently, the algebra
$\End(A_{\Phi})$ consists of
the matrix pairs
\[
[f(\Phi),f(\Phi^{*})],\qquad
f(x)\in K[x],
\]
with involution
\[
[f(\Phi),f(\Phi^{*})]^{\circ}=
[\bar f(\Phi),f(\Phi)^*].
\]
The field $T(A_{\Phi})$ can
be identified with the field
$$ K[\omega] =  K[x]/p_{\Phi}(x) K[x]$$ with involution
$f(\omega)^{\circ} = \bar
f(\omega)$. The set of
representations
$$A_{\Phi}^{f(\omega)} =
(\Phi,\Phi_1f(\Phi)),$$ where
$$0\ne f(x)=\bar
f(x)\in K[x]$$ and
$$\deg(f(x))<\deg(p_{\Phi}(x)),$$
constitute the orbit of the
representation $A_{\Phi}$.

$4^{\circ}$. From
$2^{\circ}$,  $3^{\circ}$,
and Theorem \ref{THEOREM 1},
the proof of Theorem
\ref{Theorem 6} now follows.
\end{proof}

\section{The matrices
$\sqrt[\displaystyle
*]{\Phi}$,
$\Phi_{\varepsilon}$, and
$\Phi_{(\varepsilon)}$}
\label{existtf}

Let $ K$ be a field of
characteristic not 2. In this
section we obtain the
existence conditions and
forms for the matrices
$\sqrt[\displaystyle
*]{\Phi}$,
$\Phi_{\varepsilon}$, and
$\Phi_{(\varepsilon)}$
defined in
\eqref{vrfmau1}--\eqref{vrfmau3}
by the equalities:
\begin{align*}
\sqrt[\displaystyle
*]{\Phi}&=(\sqrt[\displaystyle
*]{\Phi})^*\Phi,
   \\ 
\Phi_{\varepsilon}
&=\Phi_{\varepsilon}^*,\quad
\Phi_{\varepsilon}\Phi=
\varepsilon
(\Phi_{\varepsilon}\Phi)^*,
  \\ 
\Phi_{(\varepsilon)}
&=\varepsilon
\Phi_{(\varepsilon)}^*=
\Phi^*\Phi_{(\varepsilon)}\Phi;
\end{align*}
($\varepsilon=\pm 1$,
$\varepsilon=1$ for
nonidentity involution on
$ K$).

In the case of nonidentity
involution on $ K$, we
choose a fixed nonzero
element
\begin{equation}\label{ser23}
 k=-\bar k\ne 0;
\end{equation}
we can take $k=a-\bar a$ with
any $a\ne\bar a\in K$.

By $\Phi$ we denote an
$n\times  n$ Frobenius block,
and by
\begin{equation}\label{ser24}
\chi(x):=p(x)^s=\alpha_0x^n+
\alpha_1x^{n-1}+\dots+\alpha_n,
\end{equation}
\begin{equation}\label{ser25}
\mu(x):=p(x)^{n-1}=\beta_0x^t+
\beta_1x^{t-1}+\dots+\beta_t
\end{equation}
($\alpha_0=\beta_0=1$) we
denote the characteristic
polynomial of $\Phi$ and its
maximal divisor.

Let $$f(x)= \gamma_0x^m +
\gamma_1x^{m-1}+\dots+\gamma_m\in
 K[x].$$ A sequence
$$(a_q, a_{q+1},\dots, a_r)$$
of elements of $ K$
will be called
\emph{$f$-recurrent} if
\[
\gamma_0a_{l+m} +
\gamma_1a_{l+m-1}+\dots+\gamma_m
a_{l}=0
\]
($q\le l \le r - m$); the
sequence is completely
determined, assuming
$\gamma_0\ne 0\ne \gamma_m$,
by any fragment of length
$m$. The sequence will be
called \emph{strictly
$\chi$-recurrent} if it is
$\chi$-recurrent but not
$\mu$-recurrent (see
(\ref{ser24}) and
(\ref{ser25})).

\begin{lemma}\label{LEMMA 9}
The following two conditions
on a matrix $A$ are
equivalent:
\begin{itemize}
  \item[\rm(a)]
$A=\Phi^*A\Phi$ and $A$ is
nonsingular.
  \item[\rm(b)]
$A = [a_{j-i}]$, where the
sequence
$(a_{1-n},\dots,a_{n-1})$ is
strictly $\chi$-recurrent,
with
$\chi(x)=\chi^{\vee}(x)$.
\end{itemize}
\end{lemma}

\begin{proof}
$(a)\Longrightarrow (b)$.
Suppose the matrix $A =
[a_{ij}]$ satisfies condition
(a). Then
$$A\Phi^{-1}A^{-1}=\Phi^*,$$
and \[
   \chi(x) = \det(xI-\Phi^{*-1})=
\det(-\Phi^{*-1})\cdot
x^n\cdot \det(x^{-1}I
-\Phi^{*})=\chi^{\vee}(x).
\]
Since $$\Phi^*A\Phi=
\Phi^*[a_{i,j+1}]= [a_{i+1,j+
1}]$$ (the entries
$a_{i,n+1}$ and $a_{n+1,j}$
are defined by this
equality), we have $a_{ij}=
a_{i+1,j+ 1}$, so that the
matrix entries depend only
on the difference of the
indices; i.e., $A =
[a_{j-i}]$. That the sequence
$(a_{1-n},\dots,a_{n-1})$ is
$\chi$-recurrent follows from
the equality $A\Phi=
[a_{j-i+1}]$. Furthermore,
the recurrence is strict;
otherwise, we should have
$$(0,\dots,0,\beta_0,\dots,\beta_t)
A= 0$$ (see (\ref{ser25})),
contradicting the assumption
that $A$ is nonsingular.
\medskip

$(a)\Longleftarrow (b)$.
Suppose (b) is satisfied.
Then $$\Phi^*A\Phi=
\Phi^*[a_{j-i+1}]=
[a_{j-i}]=A.$$ We show now
that $A$ is nonsingular.
Suppose that, on the
contrary, its rows
\[
v\Phi^{n-1},\
v\Phi^{n-2},\dots,
v,\quad\text{where }v=
(a_{1-n},\dots, a_0),
\]
are linearly dependent. Then
$vf(\Phi) = 0$ for some
polynomial $f(x)\ne 0$ of
degree less than $n$. Since
$v\chi(\Phi) = 0$, we have
$vp(\Phi)^r = 0$, where
$p(x)^r$ is the greatest
common divisor of the
polynomials $f(x)$ and
$\chi(x)$. But then
\[
v\Phi^i\mu(\Phi)=
(0,\dots,0,\beta_0,\dots,\beta_t,
0,\dots,0)A=0
\]
($0\le i < n-t$; see
(\ref{ser25})); so the
sequence
$(a_{1-n},\dots,a_{n-1})$ is
$\mu$-recurrent,
contradicting condition (b).
\end{proof}

\begin{theorem}\label{THEOREM 7}
Existence conditions for the
$n\times n$ matrix
$\sqrt[\displaystyle
*]{\Phi}$ are:
\begin{itemize}
  \item[\rm(Al)]
$\chi(x) = \chi^{\vee}(x)$.

  \item[\rm(A2)]
$p(x)\ne x + (-1)^{n-1}$ in
the case of the identity
involution.
\end{itemize}
With these conditions
satisfied, we can take
\[
\sqrt[\displaystyle *]{\Phi}=
[a_{j-i}],
\]
where the sequence
$(a_{1-n},\dots,a_{n-1})$ is
$\chi$-recurrent, and is
defined by the fragment
\begin{equation}\label{ser26}
 (a_{-m},\dots,a_{m-1})
 =(\bar a,0,\dots,0,a)
\end{equation}
of length either $n$ or $n+
1$, in which
\begin{itemize}
  \item[\rm(a)]
$a = 1$ if $n = 2m$, except
for the case $p(x) =x+\alpha$
with $\alpha^{n-1}=-1$;

  \item[\rm(b)]
$a = k$ $($see
$(\ref{ser23}))$ if $n = 2m$,
$p(x) = x+\alpha$,
$\alpha^{n-1} = -1$, and also
if $n = 2m -1$, $p(x) = x +
1$;

  \item[\rm(c)]
$a=\chi(-1)$ if $n=2m-1$,
$p(x)\ne x+1$.
\end{itemize}
\end{theorem}

\begin{proof}
$1^{\circ}$. If the matrix $A
=\sqrt[\displaystyle
*]{\Phi}$ exists, then
conditions (Al) and (A2) must
be satisfied. Indeed, in view
of the relations
\[
A = A^*\Phi =\Phi^*A\Phi
\]
(see \eqref{vrfmau1}) and
Lemma \ref{LEMMA 9},
condition (Al) is satisfied,
and the entries of the matrix
\[
[a_{j-i}]= A= A^*\Phi = [\bar
a_{i-j-1}]
\]
form a strictly
$\chi$-recurrent sequence
\begin{equation}\label{ser27}
(a_{1-n},\dots,a_{n-1})=
(\bar a_{n-2},\dots,\bar a_0,
a_0,\dots, a_{n-1}).
\end{equation}
This sequence is completely
determined by the fragment
\begin{equation}\label{ser28}
(\bar a_{m-1},\dots,\bar a_0,
a_0,\dots, a_{m-1})
\end{equation}
of length $2m$, equal either
to $n$ or to $n + 1$.

Now suppose condition (A2) is
not satisfied; i.e., that the
involution is the identity and
$$p(x) = x + (-1)^{n-1}.$$
Then the vector (\ref{ser28})
is $\mu(x)
=(x+(-1)^{n-1})^{n-1}$-recurrent.
For $n = 2m$ this is obvious;
and for $n = 2m - 1$ it
follows from the property
\[
\alpha_i=\beta_i+\beta_{i-1}=
\beta_i+\beta_{n-i},\qquad
0<i<n,
\]
of the binomial coefficients
$\alpha_i$ and $\beta_i$ (see
(\ref{ser24}) and
(\ref{ser25})), since
\begin{multline}
 2[\beta_0a_{m-1}+
 \beta_1a_{m-2}+\dots
 +\beta_{n-2}a_{m-3}
 +\beta_{n-1}a_{m-2}]\\
=(\beta_0+0)a_{m-1}+
 (\beta_1+\beta_{n-1})a_{m-2}
 +
 (\beta_2+\beta_{n-2})a_{m-3}\\
 +\dots
 +(\beta_{n-1}+\beta_1)a_{m-2}
 +(0+\beta_0)a_{m-1}\\
 =\alpha_0a_{m-1}+
 \alpha_1a_{m-2}+\dots+
 \alpha_na_{m-1}=0 \label{ser29}
\end{multline}
in view of the
$\chi$-recurrence of
(\ref{ser28}). But then its
$\mu$-recurrent extension
coincides with (\ref{ser27}),
contradicting the strict
$\chi$-recurrence of
(\ref{ser27}).
\medskip

$2^{\circ}$.  If conditions
(Al) and (A2) are satisfied,
then the matrix
$\sqrt[\displaystyle
*]{\Phi}$ exists. Indeed, let
us verify that the vector
(\ref{ser26}) is strictly
$\chi$-recurrent.

\begin{itemize}
  \item
Suppose $n = 2m$. Since
(\ref{ser26}) is of length
$n$, it suffices to verify
that it is not
$\mu$-recurrent. If
$\deg(\mu(x))<n -1$, this is
obvious. If $\deg(\mu(x))=n
-1$, then the polynomial
$\mu(x)$ is of the form
$(x+\alpha)^{n-1}$, and
therefore
\[
a + \beta_{n-1}\bar a = a +
\alpha^{n-1}\bar a\ne 0.
\]

  \item
Suppose $n = 2m - 1$. Since
(\ref{ser26}) is of length $n
+ 1$, it suffices to verify
that it is $\chi$-recurrent,
i.e., that $a +\alpha_{n}\bar
a= 0$ (see (\ref{ser24})).
Condition (Al) implies that
$\alpha_n
=\bar\alpha_n^{-1}$, and so,
since $$\chi^{\vee}(x) =
\bar\alpha_n^{-1}x^n\bar
\chi(x^{-1}),$$ that
$$\chi(-1) =
-\alpha_n\overline{
\chi({-1})}.$$ If $\chi(-1)=
0$, then
$$\chi(x)=(x+1)^n,\qquad
k+\alpha_n\bar k=0.$$
\end{itemize}
Thus, the vector
(\ref{ser26}) is strictly
$\chi$-recurrent, and its
$\chi$-recurrent extension
has, in view of (Al), the
form (\ref{ser27}).
Consequently, $$A= [a_{j-i}]
= A^*\Phi.$$ By Lemma
\ref{LEMMA 9}, the matrix $A$
is nonsingular, and it can be
taken to be
$\sqrt[\displaystyle
*]{\Phi}$.
\end{proof}

\begin{theorem}\label{THEOREM 8}
Existence conditions for the
$n\times n$ matrix
${\Phi}_{\varepsilon}$ are:
\begin{itemize}
  \item[\rm(Bl)]
$\chi(x) =
\varepsilon^n\bar\chi
(\varepsilon x)$.

  \item[\rm(B2)]
$\chi(x)\notin\{x^2,x^4,x^6,\ldots
\}$ if $\varepsilon=-1$.
\end{itemize}
With these conditions
satisfied, we can take
\[
{\Phi}_{\varepsilon}=
[\varepsilon^i a_{i+j}],
\]
where the sequence
$(a_2,a_3,\dots,a_{2n})$ is
$\chi$-recurrent, and is
defined by the fragment
\begin{equation}\label{msu}
(a_2,\dots,a_{n+1})=
  \begin{cases}
    (1,0,\dots,0) & \text{if $\Phi$
    is nonsingular}, \\
    (0,\dots,0,1) & \text{if $\Phi$
    is singular}.
  \end{cases}
\end{equation}
\end{theorem}

\begin{proof}
$1^{\circ}$. Suppose
${\Phi}_{\varepsilon}$
exists. Then $$\Phi=
\Phi_{\varepsilon}^{-1}
({\varepsilon}\Phi^*)
\Phi_{\varepsilon}$$ (see
\eqref{vrfmau2}), and this
gives condition (Bl):
\[
\chi(x) =\deg(xI-\varepsilon
\Phi^*)=
\varepsilon^n\bar\chi
(\varepsilon x).
\]
Condition (B2) follows from
\eqref{vrfmau2} and Lemma
\ref{LEMMA 8}. \medskip

$2^{\circ}$.  Suppose
conditions (B1) and (B2) are
satisfied. The matrix
${\Phi}_{\varepsilon}=
[\varepsilon^i a_{i+j}]$,
defined in the statement of
Theorem \ref{THEOREM 8}, is
nonsingular. Let us verify
that it satisfies
(\ref{vrfmau2}).

If $\Phi$ is singular, this
is obvious. Suppose $\Phi$ is
nonsingular. Then the
$\chi$-recurrence of the
sequence $(a_2, a_3,\dots,
a_{2n})$ implies that
$$\Phi_{\varepsilon}\Phi =
[\varepsilon^ia_{i+j+1}];$$
and so relations
(\ref{vrfmau2}) can be
written in the form
\[
\varepsilon^ia_{i+j}=
\varepsilon^j\bar
a_{j+i},\qquad
\varepsilon^ia_{i+j+1}=
\varepsilon\varepsilon^j\bar
a_{j+i+1},
\]
i.e.,
\begin{equation}\label{ser30}
a_t= \varepsilon^t\bar
a_t,\qquad 2\le t\le 2n.
\end{equation}
We argue now by induction.
Relation (\ref{ser30})
certainly holds for $t \le n
+ 1$  (see \eqref{msu}).
Assuming it holds for $t < n
+ l$ $(l\ge 2)$, we must
verify it for $t = n + l$.
And indeed, using the
$\chi$-recurrence of the
sequence $(a_2,\dots,
a_{2n})$ and equalities
(\ref{ser24}) and (B1), we
find that
\begin{align*}
a_{n+l}&=-\alpha_1a_{n+l-1}-\dots
-\alpha_na_l\\
&=-\varepsilon\bar\alpha_1
\varepsilon^{n+l-1}\bar
a_{n+l-1}
-\dots-\varepsilon^n\bar\alpha_n
\varepsilon^l\bar
a_{l}=\varepsilon^{n+l}\bar
a_{n+l}.
\end{align*}
\end{proof}

\begin{theorem}\label{THEOREM 9}
Existence conditions for the
$n\times n$ matrix
$\Phi_{(\varepsilon)}$ are:
\begin{itemize}
  \item[\rm(Cl)]
$\chi(x) = \chi^{\vee}(x)$.

  \item[\rm(C2)]
If the involution on $ K$ is the identity and
$\varepsilon=(-1)^n$, then
$\deg(p(x)) > 1$ $($see
\eqref{ser24}$)$.
\end{itemize}
With these conditions
satisfied, we can take
\[
\Phi_{(\varepsilon)}=
[a_{j-i}],
\]
where the sequence
$(a_{1-n},\dots,a_{n-1})$ is
$\chi$-recurrent, and is
defined by the fragment
 $v=(a_{-m},\dots,a_{m})$
of length either $n$ or $n+
1$, that equals to
\begin{itemize}
  \item[\rm(a)]
$(\varepsilon\bar{\alpha}_n-1,0,
\dots,0,\alpha_n-\varepsilon)$
if $n=2m$,
$\alpha_n\ne\varepsilon$
$($see \eqref{ser24}$)$;

  \item[\rm(b)]
$(\alpha_1,-1,0,
\dots,0,-1,\alpha_1)$
$(v=(\alpha_1,-2,\alpha_1)$
for $n=2)$ if $n=2m$,
$\varepsilon=1$, and the
involution on $ K$ is
the identity;

  \item[\rm(c)]
$(-k, 0,\dots , 0, k)$ $($see
\eqref{ser23}$)$ if $n = 2m$,
$\alpha_n = 1$, and the
involution is nonidentity,
and also if $n = 2m + 1$,
$p(x) = x + \alpha$,
$\alpha^{n-1}=-1$;

  \item[\rm(d)]
$(\varepsilon, 0,\dots , 0,
1)$ if $n = 2m + 1$, in any
other case besides $p(x) = x
+ \alpha$, $\alpha^{n-1}=-1$.
\end{itemize}
\end{theorem}

\begin{proof}
$1^{\circ}$. If the matrix $A
=\Phi_{(\varepsilon)}$
exists, then conditions (Cl)
and (C2) are satisfied.
Indeed, in view of the
relations \eqref{vrfmau3} and
Lemma \ref{LEMMA 9},
condition (Cl) is satisfied,
and the entries of the matrix
$$A = [a_{j-i}] = \varepsilon
A^*$$ form a strictly
$\chi$-recurrent sequence
\begin{equation}\label{ser31}
(a_{1-n},\dots,a_{n-1})=
(\varepsilon\bar
a_{n-1},\dots,\varepsilon\bar
a_0=a_0,\dots, a_{n-1})
\end{equation}

Suppose condition (C2) is not
satisfied. By (Cl), $$p(x) =
p^{\vee}(x) = x \pm 1,$$ and
the fragment $(\varepsilon
a_m,\dots,a_{m})$ of length
either $n$ or $n + 1$, of the
vector (\ref{ser31}) is
$\mu$-recurrent. This is
obvious if $n = 2m + 1$ since
$\varepsilon = -1$; and if $n
= 2m$, it follows from
(\ref{ser29}) as applied to
the fragment (replace $m$ in
(\ref{ser29}) by $m + 1$).
But then the vector
(\ref{ser31}) is also
$\mu$-recurrent, and we have
a contradiction.

$2^{\circ}$. If conditions
(Cl) and (C2) are satisfied,
then $\Phi_{(\varepsilon)}$
exists. To show this, let us
verify that the vector $v$ of
Theorem \ref{THEOREM 9} is
strictly $\chi$-recurrent and
of the form
\[
(\varepsilon\bar
a_m,\dots,\varepsilon\bar
a_0=a_0,\dots, a_m).
\]
\begin{itemize}
  \item
The vector in (a) is
$\chi$-recurrent, since its
length is $n + 1$ and, by
(Cl),
$\alpha_n\bar{\alpha}_n= 1$.

  \item
The vector in (b) is
$\chi$-recurrent, since for
the identity involution
conditions (Cl) and (C2)
imply
\[
\chi(1)
=\alpha_n^{-1}\chi(1)\ne 0,
\quad \alpha_n=1,\quad
\alpha_{n-1}=\alpha_1.
\]
The vector is not
$\mu$-recurrent, since $t\le
n - 2$ (by (\ref{ser25}) and
(C2)) and $\beta_t= 1$ (by
the equality $p(x) =
p^{\vee}(x)$ and (C2)).

  \item
If $n = 2m + 1$, $p(x) = x +
\alpha$, and $\alpha^{n-1}=
-1$ (see (c)), then the
involution is nonidentity:
otherwise $$p(x) =
p^{\vee}(x) = x \pm 1,$$
contradicting the equality
$\alpha^{n-1}= -1$.

  \item
The vector in (d) is not
$\mu$-recurrent, in view of
(C2).
\end{itemize}

Now let (\ref{ser31}) be the
$\chi$-recurrent extension of
the vector $v$. Then the
matrix $A =[a_{j-i}]$ is
equal to $\varepsilon A^*$,
and by Lemma \ref{LEMMA 9} it
can be taken for
$\Phi_{(\varepsilon)}$.
\end{proof}

\end{document}